\documentclass[]{article}

\usepackage{arxiv}

\usepackage{amssymb}
\usepackage{amsfonts}
\usepackage{amsmath,amsthm}
\usepackage{stmaryrd}
\usepackage{rotating}
\usepackage{graphicx,color}
\usepackage{subcaption}
\usepackage[dvipsnames]{xcolor}
\usepackage{epstopdf}
\usepackage{tikz}
\usepackage{array}
\usepackage[sans]{dsfont}
\usepackage{booktabs}
\usepackage{multicol}
\usepackage{multirow}
\usepackage{gensymb}
\usepackage{lineno,hyperref}
\usepackage{arydshln}
\usepackage[export]{adjustbox}
\usepackage{caption}

\theoremstyle{remark} 
\newtheorem{remark}{\bf Remark}
\newtheorem{testcase}{\bf Test Case}

\usepackage{pgfplots}
\pgfplotsset{compat = 1.10}
\definecolor{mycolor1}{rgb}{1.00000,1.00000,0.00000}%
\usetikzlibrary{decorations.markings}
\usetikzlibrary{shapes,arrows}
\usetikzlibrary{shapes.geometric}
\usetikzlibrary{fit}					
\usetikzlibrary{backgrounds}
\usetikzlibrary{patterns}

\renewcommand{\vec}[1]{\boldsymbol#1} 
\newcommand{\tensorTwo}[1]{\boldsymbol#1}

\renewcommand{\Vec}[1]{\mathsf{#1}}
\newcommand{\Mat}[1]{\mathsf{#1}}
\newcommand{\blkVec}[1]{\boldsymbol{\mathsf{#1}}}
\newcommand{\blkMat}[1]{\boldsymbol{\mathsf{#1}}}

\newcommand{\sub}[1]{_\mathit{#1}}

\usepackage{todonotes}
\pgfplotsset{
	colormap={parula}{
		rgb255=(53,42,135)
		rgb255=(15,92,221)
		rgb255=(18,125,216)
		rgb255=(7,156,207)
		rgb255=(21,177,180)
		rgb255=(89,189,140)
		rgb255=(165,190,107)
		rgb255=(225,185,82)
		rgb255=(252,206,46)
		rgb255=(249,251,14)
	},
}

\pgfplotsset{
	colormap={hotDesaturated}{
		rgb255=(71,71,219)
		rgb255=(0,0,91)
		rgb255=(0,255,255)
		rgb255=(0,127,0)
		rgb255=(255,255,0)
		rgb255=(255,96,0)
		rgb255=(107,0,0)
		rgb255=(224,76,76)
	},
}

\begin{document}
	
\title{Simulation of coupled multiphase flow and geomechanics in porous media with embedded discrete fractures}

\author{
	Matteo Cusini \\
	Atmospheric, Earth and Energy Division\\
	Lawrence Livermore National Laboratory \\
	7000 East Ave., Livermore, CA 94550, USA\\
	\texttt{cusini1@llnl.gov} \\
	\And
	Joshua A.~White \\
	Atmospheric, Earth and Energy Division\\
	Lawrence Livermore National Laboratory \\
	7000 East Ave., Livermore, CA 94550, USA\\
	\texttt{white230@llnl.gov} \\
	\And
	Nicola Castelletto \\
	Atmospheric, Earth and Energy Division\\
	Lawrence Livermore National Laboratory \\
	7000 East Ave., Livermore, CA 94550, USA\\
	\texttt{castelletto1@llnl.gov} \\
	\And
	Randolph R.~Settgast \\
	Atmospheric, Earth and Energy Division\\
	Lawrence Livermore National Laboratory \\
	7000 East Ave., Livermore, CA 94550, USA\\
	\texttt{settgast1@llnl.gov} \\
}

\maketitle

\begin{abstract}
	In fractured natural formations, the equations governing fluid flow and geomechanics are strongly coupled. Hydrodynamical properties depend on the mechanical configuration, and they are therefore difficult to accurately resolve using uncoupled methods. In recent years, significant research has focused on discretization strategies for these coupled systems, particularly in the presence of complicated fracture network geometries. In this work, we explore a finite-volume discretization for the multiphase flow equations coupled with a finite-element scheme for the mechanical equations. Fractures are treated as lower dimensional surfaces embedded in a background grid.  Interactions are captured using the Embedded Discrete Fracture Model (EDFM) and the Embedded Finite Element Method (EFEM) for the flow and the mechanics, respectively. This non-conforming approach significantly alleviates meshing challenges. EDFM considers fractures as lower dimension finiten volumes which exchange fluxes with the rock matrix cells. The EFEM method provides, instead, a local enrichment of the finite-element space inside each matrix cell cut by a fracture element. Both the use of piecewise constant and piecewise linear enrichments are investigated.  They are also compared to an Extended Finite Element (XFEM) approach. One key advantage of EFEM is the element-based nature of the enrichment, which reduces the geometric complexity of the implementation and leads to linear systems with advantageous properties. Synthetic numerical tests are presented to study the convergence and accuracy of the proposed method. It is also applied to a realistic scenario, involving a heterogeneous reservoir with a complex fracture distribution, to demonstrate its relevance for field applications.
\end{abstract}


\setcounter{footnote}{0} 

\section{Introduction}
In many geoengineering applications, the decision-making process is supported by numerical simulations, e.g. oil and gas fields, geothermal plants, or carbon storage reservoirs. Such systems can only be operated efficiently and safely with a thorough understanding of the flow and transport processes in the subsurface and how they interact with the mechanical response. Accurate numerical models are an important tool to evaluate the risk of undesirable phenomena such as early breakthrough, poor sweep efficiency, excessive subsidence, loss of containment, or induced seismicity.

Many reservoirs are either naturally or deliberately fractured. In fractured formations, the coupling between fluid flow and geomechanics is particularly strong as the hydrodynamical properties (e.g. permeability, storage) of the fractures are strongly dependent on the mechanical configuration.  This makes numerical simulations of such systems particularly challenging. From a mathematical point of view, one has to solve a system of coupled partial differential equations. In particular, the flow of fluids and the mechanical response are described by mass and momentum balances, along with several nonlinear constitutive relationships \cite{Wang2000}. Numerical difficulties frequently arise from the geometrical complexity of real fracture networks.    

There exist two broad classes of discretization methods which model fractures as lower dimensional entities (e.g., 2D surfaces in a 3D domain): conforming-grid methods and non-conforming (or embedded) methods. The first class relies on a single grid that follows the geometry of the fractures so that conventional discretization approaches may be applied.  This includes finite-element (FE), finite volume (FV), and combined FE/FV schemes \cite{karimi2003DFN,Garipov2016,Glaser2017,Randy2017,Brenner2020}. Embedded methods, instead, consider the fractures as independent surfaces overlain on a separate background grid. This approach circumvents the complexity of generating a single grid that honors the geometry of the fracture network. Separate, non-conforming grids are generated for the rock matrix and for each fracture. Such approaches have been proposed for both FV discretization of flow \cite{Lee2001,Hajibeygi2011,MousaFADM2018} and poromechanics \cite{Deb2017,Ucar2018,Deb2020} and as extended finite-element methods (XFEM) for flow \cite{AlessioAnnaXFEM2013,Berrone2014}, mechanics, and poromechanics \cite{Rethore2007,Berrone2014,Khoei2014,Giovanardi2017}. Recently, these methods have also been combined for simulation of coupled multiphase flow and mechanics in fractured porous media using a mixed FE/FV discretization \cite{Ren2016}. 

In this work we propose a FV discretization of the multiphase mass balance equations along with a FE scheme for the momentum balance equation. The contribution of the fractures to the flow and transport is captured by employing the embedded discrete fracture model (EDFM) \cite{Lee2001,Hajibeygi2011} which is an extension of a traditional two-point flux finite-volume scheme to discontinuous media. The contact/separation problem at each fracture element is handled using the Embedded Finite Element (EFEM) method.  This method employs an element-local enrichment of the FE space using the concept of Assumed Enhanced Strain \cite{Simo1990,FosterAES,Wells2001,Oliver2006,Borja2008,Wu2011,Deb2017}.  Both  piecewise constant\cite{Simo1990} and linear \cite{Linder2007} interpolation of the displacement jump within an element are considered. While the focus in this paper is the EFEM framework, we also describe a comparable XFEM-based method to highlight the overlap and differences in the resulting approach.  EFEM and XFEM are often viewed as competing methods, though they share many similarities \cite{Borja2008}.  In fact, a well-written code can implement both in a straightforward manner with significant code re-use.  

The paper is organized as follows. The governing equations and constitutive assumptions are presented in \S\ref{Sec:ProblemStatement}. In \S\ref{Sec:Discretization}, we describe the combined FE/FV discretization, with a focus on the EFEM and EDFM enrichment strategy. The accuracy and consistency of the method are studied in \S\ref{Sec:NumResults} through simple benchmarks. More complex numerical experiments are then conducted to demonstrate the applicability of the method to realistic field applications. Finally, concluding remarks are provided in \S\ref{Sec:Conclusions}.

\section{Problem Statement}
\label{Sec:ProblemStatement}

The goal of this work is to model tight coupling between multiphase flow and elastic deformation in a porous and fractured medium.  For simplicity, we present the formulation and discretization for the two-dimensional case, but the extension to three-dimensions is reasonably straightforward. We assume quasi-static, small-strain kinematics throughout.

\begin{figure}[t]
	\hfill
	\begin{subfigure}[c]{.35\linewidth}
		\centering
		\includegraphics[width=\linewidth]{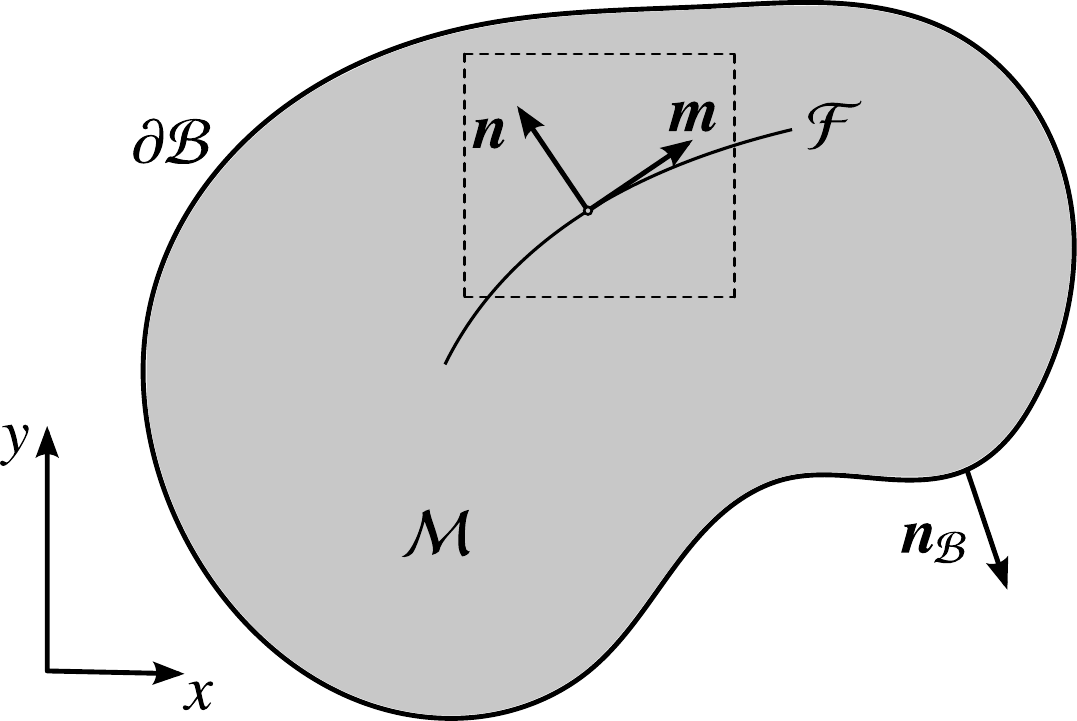}
		\caption{}
	\end{subfigure}
	\hfill
	\begin{subfigure}[c]{.35\linewidth}
		\centering
		\includegraphics[width=\linewidth]{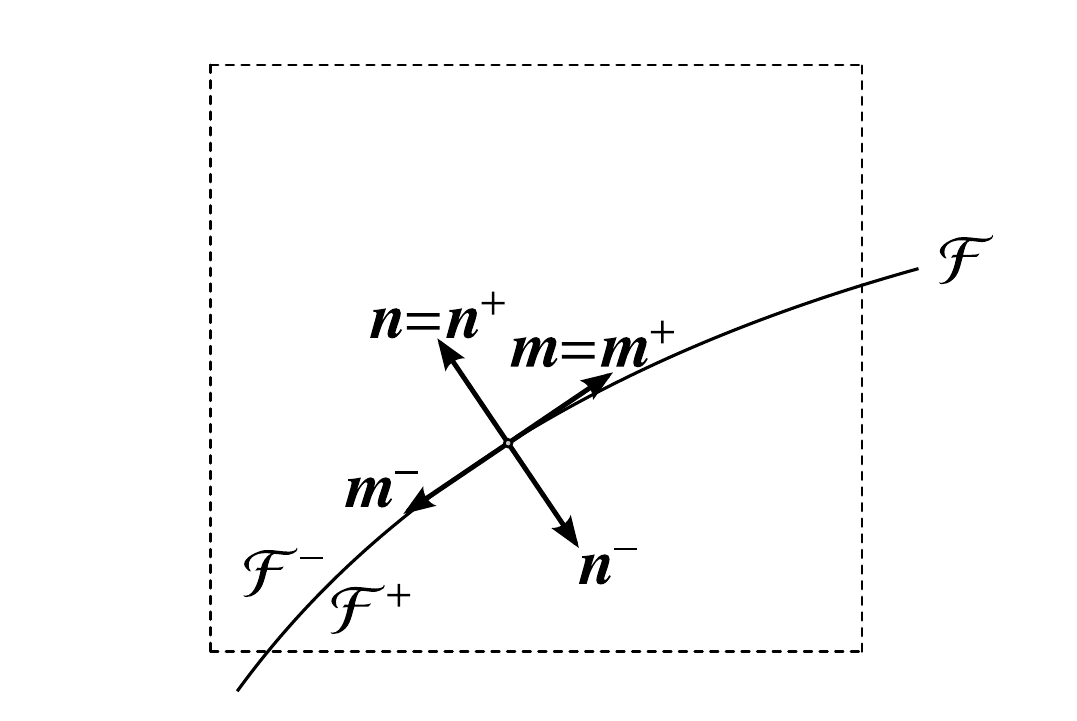}
		\caption{}
	\end{subfigure}
	\hfill\null
	\caption{Problem geometry.}
	\label{Fig:Geometry}
\end{figure}

Let $\mathcal{B} \subset \mathbb{R}^2$ be an open domain (Figure~\ref{Fig:Geometry}).  Its external boundary is  $\partial \mathcal{B}$, with outward normal vector $\vec{n}_\mathcal{B}$.  For the application of boundary conditions, the external boundary is divided into non-overlapping portions $\partial \mathcal{B} = \overline{\partial \mathcal{B}_u} \,\cup\,\overline{\partial \mathcal{B}_t} = \overline{\partial \mathcal{B}_p} \,\cup\,\overline{\partial \mathcal{B}_q}$ where Dirichlet and Neumann conditions for the mechanical and flow portions of the problem will be applied.  The continuous body is internally cut by one or more fractures which form an embedded, lower-dimensional domain $\mathcal{F}$ with boundary points $\partial \mathcal{F}$.  For simplicity, it is assumed fractures do not intersect the external boundary.  We will commonly refer to $\mathcal{F}$ as the fracture network, anticipating the general case of multiple, interconnected fracture segments.  We assume also that the network geometry can be well approximated by one or more polylines.  Each fracture segment has two faces, with the set of positive and negative faces denoted as $\mathcal{F}^\pm$.  By convention, the local normal is chosen as $\vec{n} = \vec{n}^+ = -\vec{n}^-$.  Similarly, we let $\vec{m}$ denote the tangent vector forming a local, right-handed coordinate system at the fracture surface.  The porous matrix is then $\mathcal{M}=\mathcal{B} \setminus \mathcal{F}$. We model system behavior in the time domain $\mathbb{T}=(0,t_\text{max}]$ from known initial conditions.  

For the mechanical deformation, the continuum matrix $\mathcal{M}$ behaves as a poroelastic medium with deformation field $\vec{u}$.   Without loss of generality, homogenous boundary conditions $\vec{u}=\vec{0}$ are prescribed on $\partial \mathcal{B}_u$.  Across $\mathcal{F}$, the deformation field is potentially discontinuous, with $ \llbracket \vec{u} \rrbracket = \vec{u}^+ - \vec{u}^-$ denoting the displacement jump.  It is convenient to partition this jump into normal and tangential components as $\llbracket \vec{u} \rrbracket = w_n \vec{n} + w_m \vec{m}$, where the scalars $w_n$ and $w_m$ are the fracture aperture and tangential slip magnitude, respectively.  

For the multiphase flow, both the porous medium and the fractures are filled with two compressible fluids, a wetting ($w$) and non-wetting ($nw$) phase.  Let $s_\pi$ and $p_\pi$ denote the phase saturation and phase pressure of fluid phase $\pi = \{ w, nw \}$.   Since the two fluids fill the voids, the saturations satisfy the closure condition $s_w + s_{nw} = 1$.   In the following, the wetting fluid phase saturation is selected as a primary unknown and will be denoted by lower case $s$ without subscript.   Capillarity effects are not considered, leading to the simplification $p_w = p_{nw} = p$.  This is a common assumption in many conventional reservoir applications, but it is not central to the method.  We emphasize that fluids may flow both in the matrix $\mathcal{M}$ and in the fracture network $\mathcal{F}$, with interchange between the two domains.

With these preliminaries, the strong form of the initial-boundary-value problem is to find the displacement $\vec{u} : \overline{\mathcal{M}} \times \mathbb{T} \rightarrow \mathbb{R}^2$, saturation $s : \overline{\mathcal{B}} \times \mathbb{T} \rightarrow \mathbb{R}$, and pressure $p : \overline{\mathcal{B}} \times \mathbb{T} \rightarrow \mathbb{R}$ such that 
\begin{subequations}\label{Eqn:StrongForm}
	\begin{align}
	&\nabla \cdot \tensorTwo{\sigma} + \rho\vec{g} = \vec{0} & &\mbox{ on } \mathcal{M} \times \mathbb{T} & &\mbox{(matrix momentum balance),} \label{Eqn:MomBal}\\
	%
	&\llbracket \tensorTwo{\sigma} \rrbracket \cdot \vec{n} = \vec{0} & &\mbox{ on } \mathcal{F} \times \mathbb{T} & &\mbox{(fracture traction balance),} \label{Eqn:TractionBal}\\
	%
	&\dot{m}^m_{\pi}+ \nabla \cdot \vec{q}^m_{\pi} - q^m_{\pi} + q^{mf}_{\pi} = 0 & &\mbox{ on } \mathcal{M} \times \mathbb{T} & &\mbox{(matrix mass balance for fluid phase $\pi=\{w,nw\}$),} \label{Eqn:MatrixMassBal}	 \\
	%
	&\dot{m}^{f}_{\pi}+ \blacktriangledown \cdot \vec{q}^f_{\pi} - q^f_{\pi} - q^{mf}_{\pi} = 0 & &\mbox{ on } \mathcal{F} \times \mathbb{T} & &\mbox{(fracture mass balance for fluid phase $\pi=\{w,nw\}$),} \label{Eqn:FracMassBal}	 \\
	%
	\intertext{subject to boundary conditions}
	&\vec{u} = {\vec{0}} & &\mbox{ on } \partial \mathcal{B}_u \times \mathbb{T} & &\mbox{(prescribed displacement),} \label{eq:momentumBalanceS_DIR}\\
	&\tensorTwo{\sigma} \cdot \vec{n}_{\mathcal{B}} = \bar{\vec{t}} & &\mbox{ on } \partial \mathcal{B}_t \times \mathbb{T} & &\mbox{(prescribed total traction),} \\
	&p =\bar{p} & &\mbox{ on }  \partial \mathcal{B}_p \times \mathbb{T} & &\mbox{(prescribed pore pressure),}\\    
	&s =\bar{s} & &\mbox{ on }  \partial \mathcal{B}_p \times \mathbb{T} & &\mbox{(prescribed wetting phase saturation),}\\        
	&\vec{q}^m_{\pi} \cdot \vec{n}_{\mathcal{B}} = \bar{q}^m_\pi & &\mbox{ on }  \partial \mathcal{B}_q \times \mathbb{T} & &\mbox{(prescribed mass flux for phase $\pi=\{w,nw\}$),}\\           
	%
	\intertext{and initial conditions}
	&\vec{u}(\vec{x}, 0) = \vec{u}_0 (\vec{x}) & &\mbox{ } \vec{x} \in \overline{\mathcal{M}} & &\mbox{(initial displacement),}\\
	&s(\vec{x}, 0) = s_0 (\vec{x}) & &\mbox{ } \vec{x} \in \overline{\mathcal{B}} & &\mbox{(initial wetting phase saturation),}\\
	&p(\vec{x}, 0) = p_0 (\vec{x}) & &\mbox{ } \vec{x} \in \overline{\mathcal{B}} & &\mbox{(initial pore pressure).}
	\end{align}
\end{subequations}
In these equations, the following variables, operators, and constitutive relationships are introduced:
\begin{itemize}
	
	\item The total Cauchy stress tensor is $\vec{\sigma} = \mathbb{C}_{dr}: \nabla^s \vec{u} - b p \vec{1}$, where $\mathbb{C}_{dr}$ is a fourth-order tensor of drained elastic moduli, $\nabla^s$ is the symmetric gradient operator, $b \in (\phi^0,1]$ is Biot's coefficient, with $\phi^0$ the reference porosity, and $\vec{1}$ is the second-order unit tensor.  For an isotropic model, $\mathbb{C}_{dr}$ can be expressed in terms of the drained skeleton modulus $K_{dr}$ and Poisson ratio $\nu$.
	
	\item The mixture density is computed as $\rho = (1-\phi) \rho_s + s \phi \rho_w + (1-s) \phi \rho_{nw}$ using matrix porosity $\phi$ and individual phase densities for the solid, wetting phase, and non-wetting phases. The gravitational vector is denoted by $\vec{g}$.
	
	\item The phase densities follow the compressible model $\rho_\pi = \rho_\pi^0 \exp [(p-p_\pi^0)/K_\pi]$ with reference density $\rho_\pi^0$ at reference pressure $p_\pi^0$ and phase bulk modulus $K_\pi$.
	
	\item Porosity changes depend on displacement and pore pressure as $\dot \phi = b\, \nabla \cdot \dot{\vec{u}} + \dot{p} / N$ with $N= K_{dr} / [(b-\phi^0)(1-b)]$.  
	
	\item In this work, we consider both open and closed fractures, and a suitable constitutive model for the fracture traction must be provided.  Equation \eqref{Eqn:TractionBal} expresses a traction balance on the fracture surface. For it to be satisfied, the total traction on the fracture must be equal to $\vec{t} =  \vec{\sigma}^+ \cdot \vec{n} = - \vec{\sigma^-} \cdot \vec{n}$.  This traction is additively decomposed as $\vec{t} = \vec{t}' - p \vec{n}$, where $\vec{t}'$ is the effective traction associated with mechanical contact, and the second term is a normal traction created by the fluid pressure in the fracture.  We assume 100\% of the fluid pressure is mobilized in creating traction on the fracture walls. When the fracture is open ($w_n > 0$) the mobilized traction depends only on the fluid pressure and $\vec{t}'=\vec{0}$.  When the fracture closes, a no-interpenetration constraint $w_n=0$ is enforced.  In the closed state, tangential tractions can also be generated, which are modeled using a regularized Coulomb model to describe frictional sliding.  The resulting nonlinear model can be expressed in a general rate form $\dot{\vec{t}'} = \mathbb{D} \cdot \llbracket \dot{\vec{u}} \rrbracket$, where $\mathbb{D}$ is a second-order tensor of tangent moduli and $ \llbracket \dot{\vec{u}} \rrbracket$ is the velocity jump. See \cite{White2014} for implementation details. 
	
	\item For modeling fluid flow in the fracture, we make a distinction between the mechanical aperture $w_n$ and the hydraulic-aperture $w_h = w_n + w_0$.  The basic concept is that when two rough surfaces are in contact, the voids between asperities provide a pathway for fluid flow even when the mechanical aperture is nominally zero.  The correction term $w_0 > 0$ allows for fluid storage and flow even under contact conditions.  For more accurate results in realistic scenarios, more complicated stress-dependent fracture closure models could also be considered \cite{Witherspoon1980,barton1985strength}.
	
	\item In the matrix, the phase mass per unit volume is $m_\pi^m = \phi s_\pi \rho_\pi$.  The porosity introduces a deformation coupling in the mass balance equations.  The phase mass flux is $\vec{q}_\pi^m = \rho_\pi \vec{v}_\pi$ with phase velocity $\vec{v}^m_\pi=- \lambda_\pi \vec{\kappa} \cdot \nabla (p + \rho_\pi g z)$ following the generalized Darcy's law \cite{MusMer36}.  The absolute permeability tensor is $\vec{\kappa}$.  The phase mobility $\lambda_\pi = k_{r\pi}/\mu_\pi$ is a function of the relative permeability relationship $k_{r\pi} (s_\pi)$ and the phase viscosity $\mu_\pi$.   The gravitational acceleration is $g$ and  $z$ is the elevation above a datum.   In the numerical examples, we employ a quadratic relative permeability model and constant viscosities.  
	
	\item In the fracture, we neglect the presence of any infilling material.  The phase mass per unit surface area is $m_\pi^f = w_h s_\pi \rho_\pi$.  The hydraulic aperture $w_h$ is used here rather than the mechanical aperture $w_n$.   Introducing the tangential projection matrix $\Pi = (\tensorTwo{1} - \vec{n} \otimes \vec{n}$), the operators $\blacktriangledown() = \Pi \cdot \nabla()$ and $\blacktriangledown \cdot () = \Pi : \nabla()$  are the \emph{tangential} gradient and divergence operators in the lower dimensional domain $\mathcal{F}$.  The phase flux is $\vec{q}_\pi^f = w_h \rho_\pi \vec{v}_\pi$ with phase velocity $\vec{v}_\pi^f = - \lambda_\pi \kappa^f \blacktriangledown (p + \rho_\pi g z)$.   That is, fluid flow within the fractures is driven by the tangential gradient of the hydraulic potential.  Here, the fracture permeability is $\kappa^f = w_h^2/12$ following the classic lubrication model \cite{Bear1972,Witherspoon1980,rutqvist1995coupled}. The same quadratic relative permeability curves, considered for the rock matrix, are assumed to be valid in the fractures.
	
	\item The terms $q_\pi^m$, $q_\pi^f$ are source terms used to model wells that inject or extract fluid from the rock matrix or fracture network, respectively.  The term $q_\pi^{mf}$ is an interchange term used to model the transfer of fluid from the matrix to the fracture, or vice versa.  In the embedded discretization method described below, it is more natural to model this interchange via a source term, rather than as a boundary condition that must be enforced between the matrix and fracture surfaces. Additionally, because of the choice of the EDFM formulation, only highly conductive fractures are considered in this work. This hypothesis could be relaxed, allowing for the presence of flow barriers, by considering  an extension of the EDFM formulation, namely p-EDFM \cite{Tene2017,Jiang2017}.
\end{itemize}

The model above has sufficient complexity to describe many realistic subsurface systems.  Of course, several alternative constitutive models could be introduced without changing the underlying nature of the governing equations, and certain assumptions could be relaxed to better describe particular applications.  Our primary goal in this work, however, is to test a particular discretization strategy to capture the complicated hydromechanical interactions that can result.

\section{Discretization}
\label{Sec:Discretization}

\begin{figure}[htbp]
	\centering
	\includegraphics[width=\linewidth]{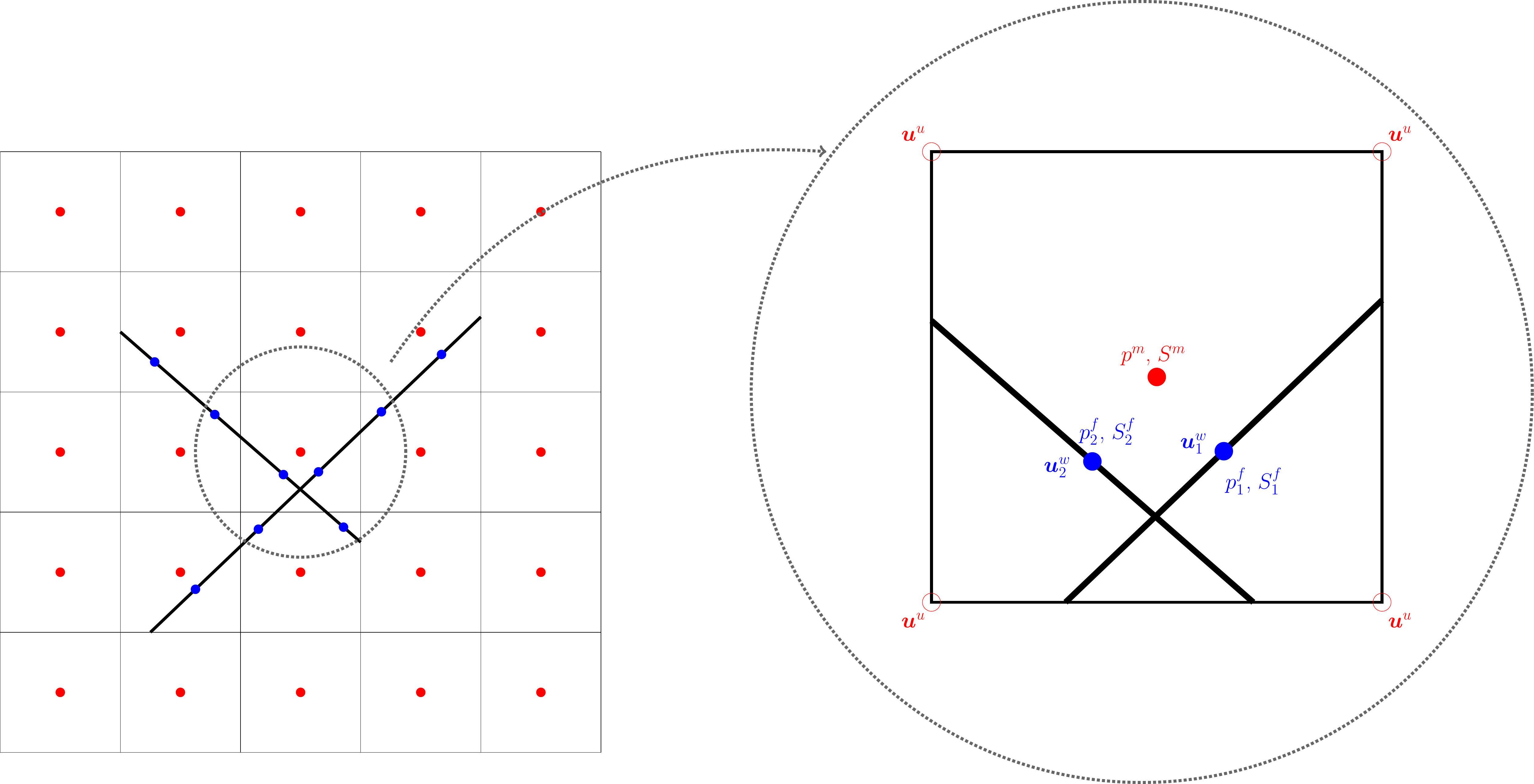}
	\caption{Grid and location of the unknowns for a two-dimensional domain.}
	\label{fig:discretization}
\end{figure}

Figure~\ref{fig:discretization} provides a simple illustration of the embedded discretization strategy adopted here.  We first partition the domain using a computational mesh $\mathcal{T}_\mathcal{B}$.  This mesh consists of non-overlapping cells $K_i$ such that $\mathcal{B} \approx  \bigcup_i K_i$.  The intersection of the fracture network with this background mesh defines a fracture triangulation $\mathcal{T}_\mathcal{F}$.  In particular, we assume $\mathcal{F} \approx  \bigcup_j k_j$, where fracture segment $k_j = \mathcal{F} \cap K_j$ for any cell $K_j$ cut by the network.  Note that we will consistently use an uppercase / lowercase notation---$(K,k)$---to indicate a cell and its corresponding fracture segment. For the moment, let us assume that each cell is cut by at most one linear segment.  The case of multiple segments intersecting a cell will be addressed once the preliminaries are established.  We assume throughout, however, that such segments completely cut the cell and do not partially penetrate.  It is convenient to denote the union of cut cells as the enriched subset of the triangulation  $\mathcal{T}_\mathcal{E} \subseteq \mathcal{T}_\mathcal{B}$. 

On a given cell, it is helpful to work in a fracture-aligned coordinate system.  To do so, we define local basis vectors $\{\vec{n},\vec{m}\}$ and origin $\vec{x}_k$, choosing $\vec{x}_k$ as the mid-point of the fracture segment. Let $\{y_n,y_m\}$ denote the normal and tangential coordinates associated with this system,
\begin{align}
y_n &= \vec{n} \cdot (\vec{x}-\vec{x}_k) \,, \notag \\
y_m &= \vec{m} \cdot (\vec{x}-\vec{x}_k)  \,.
\end{align}

Using these triangulations, the unknown fields are approximated with discrete counterparts---i.e. $\vec{u}^h$, $p^h$, and $s^h$.  To avoid a proliferation of superscripts, we will drop the standard ``h" notation and simply remark that all fields beyond this point should be understood as discrete approximations. In cells that are not cut by the fracture network, standard interpolation strategies may be adopted. Any cell cut by the fracture network, however, must be enriched with additional degrees of freedom to capture discontinuities.  

The governing equations are time-dependent and are discretized into discrete timesteps $\{0, t_1, t_2, ..., t_\text{max}\}$.  Let $\Delta t = t_\tau - t_{\tau-1}$ denote the current time interval at timestep $\tau$, and more generally $\Delta (\cdot) = (\cdot)_\tau - (\cdot)_{\tau-1}$ the discrete increment of a given quantity. The governing equations are discretized using a fully-implicit strategy, in which all unknowns are simultaneously updated as part of a Newton search.  That is, given the previous timestep solution $\{\vec{u},p,s\}_{\tau-1}$, we seek the next timestep solution $\{\vec{u},p,s\}_{\tau}$ in a monolithic fashion.  For presentation purposes, however, it is convenient to group the governing equations into a mechanics subproblem and a fluid flow subproblem.  The former relies on a finite element discretization, while the latter relies on a finite volume discretization. 

\subsection{Mechanics Discretization}

We begin by introducing two discrete spaces.  The first is the continuous bilinear finite element space,
\begin{align}
\boldsymbol{U} & := \left\{ \vec{\eta} \left| \right.  \vec{\eta} \in [C^0(\overline{\mathcal{B}})]^2,\,   \vec{\eta}_{\left| K \right.} \in [{\mathbb{Q}}_1(K)]^2 \; \forall K \in \mathcal{T}_\mathcal{B},\,
\vec{\eta} = {\vec{0}} \text{ on }\partial \mathcal{B}_u
\right\} \,.
\end{align}
Here, $C^0(\overline{\mathcal{B}})$ is the space of continuous functions on the closed  domain $\overline{\mathcal{B}}$, and $\mathbb{Q}_1(K)$ is the space of bilinear polynomials on $K$.   The space $\boldsymbol{U}$ satisfies Dirichlet boundary conditions on the displacement field.  Without loss of generality, homogeneous conditions have been assumed.  For any element $K$, let $\{N_1,N_2,N_3,N_4\}$ denote the standard ``hat'' shape functions associated to the four nodes of a quadrilateral, forming a basis for $\mathbb{Q}_1(K)$.  Any element $\vec{u} \in \boldsymbol{U}$ is locally interpolated on cell $K$ as,
\begin{equation}
\vec{u}_{\left| K \right.} = \sum_{a=1}^{4} N_a(\vec{x})\, \vec{u}_a = \sum_{b=1}^{8} u_b \,\vec{\eta}_b(\vec{x})   \,.
\end{equation}
Here, we have introduced two equivalent representations: one using four scalar shape functions $\{N_a\}$ and vector-valued weighting coefficient $\{\vec{u}_a\}$, and a second using eight vector-valued shape functions $\{\vec{\eta}_a\}$ and scalar weights $\{u_a\}$.  Note that $\vec{\eta}_1 = (N_1, 0)^T$,  $\vec{\eta}_2 = (0,N_1)^T$, ... , $\vec{\eta}_8 = (0, N_4)^T$, and the two forms are readily interchangeable.    This interpolation is used to approximate the continuous portion of the displacement field everywhere in $\overline{\mathcal{B}}$.  

We also  introduce a second space $\boldsymbol{W}$ consisting of local enrichments added to cut elements in $\mathcal{T}_\mathcal{E}$,
\begin{align}
\boldsymbol{W} & := \left\{ \vec{\phi} \left| \right.  \vec{\phi} \in [L^2(\mathcal{B})]^2,\, \vec{\phi}_{\left| K \right.} \in \boldsymbol{\mathbb{E}}(K) \; \forall K \in \mathcal{T}_\mathcal{E}, \,  \vec{\phi}_{\left| K \right.} = \vec{0} \; \forall K \notin \mathcal{T}_\mathcal{E}  
\right\} \,.
\end{align}
Here, $L^2(\mathcal{B})$ is the space of square Lebesgue-integrable functions on ${\mathcal{B}}$, 
and $\boldsymbol{\mathbb{E}}$ is a space of vectorial enrichment functions on cut elements.  These enrichments will be defined by construction.  To begin, let $\vec{\xi}(\vec{x})$ denote a continuous displacement field on $K$, with the expansion
\begin{equation}
\vec{\xi} = \sum_{b} {w}_b \, \vec{\xi}_b(\vec{x}) \,,
\end{equation}
in terms of a set of basis vectors $\{\vec{\xi}_b\}$ and enrichment weights $\{w_b\}$, as yet undefined.  Locally, the enriched displacement field is approximated as,
\begin{align}
\vec{u}_{\left| K \right.} & = \widetilde{\vec{u}}+ H \, \vec{\xi} \,
\label{eqn:hdecomp}
\end{align}
where $\widetilde{\vec{u}}$ is a locally-continuous displacement field, and $H(\vec{x})$ is a Heaviside function centered on the fracture segment $k$.  We observe that $H \, \vec{\xi}$ is an additional displacement component that is added to the positive side of the fractured element, inducing a displacement jump across $k$ and allowing for separate motion of the two sides of the element.  A variety of enrichments  $\{\vec{\xi}_b\}$ could be imagined, but here we explore three specific bases, which we will denote as EFEM(0), EFEM(1), and XFEM.  

\paragraph{EFEM(0)}
The simplest basis consists of piecewise constant enrichments for the normal and tangential displacement,
\begin{align}
\vec{\xi}_1 = \vec{n} \,,
\qquad \vec{\xi}_2 = \vec{m} \,.
\label{eq:constbasis}
\end{align}
This is the most common basis used in EFEM methods, going back to the earliest formulations \cite{Simo1990}.  It allows for piecewise constant opening and slip of the fracture segment, but no rotations or relative stretching.

\paragraph{EFEM(1)}
The second basis allows for linear displacement jumps, using the specific enrichments,
\begin{align}
\vec{\xi}_1 = \vec{n} \,,
\qquad \vec{\xi}_2 = \vec{m} \,,
\qquad \vec{\xi}_3 = y_m\vec{m} \,,
\qquad \vec{\xi}_4 = y_m\vec{n}-y_n\vec{m} \,,
\label{eq:linderbasis}
\end{align}
This particular basis was proposed in \cite{Linder2007} and further developed in subsequent works \cite{Linder2009,Armero2009,Linder2013}.  The first two modes represent rigid body motions in the normal and tangential directions.  The third provides a tangential stretching mode, while the fourth is a rigid rotation.  Due to the spatial variation of the enrichment, the displacement jump may vary linearly along the fracture length, providing a better approximation of the jump kinematics. 

\paragraph{XFEM}
A limitation of the previous enrichments is that they only allow for a subset of the separation and deformation modes possible for a separating bilinear element.  While the dominant modes are captured, the two-sides of the fracture are not perfectly independent.  Full separation can be achieved by introducing a complete, eight-mode basis,
\begin{equation}
\renewcommand{\arraycolsep}{20pt}
\begin{matrix}
\vec{\xi}_{1} = N_1 \, \vec{n} \,, 
&\vec{\xi}_{2} = N_2 \, \vec{n} \,, 
&\vec{\xi}_{3} = N_3 \, \vec{n} \,, 
&\vec{\xi}_{4} = N_4 \, \vec{n} \,,  \\
\vec{\xi}_{5} = N_1 \, \vec{m} \,, 
&\vec{\xi}_{6} = N_2 \, \vec{m} \,, 
&\vec{\xi}_{7} = N_3 \, \vec{m} \,, 
&\vec{\xi}_{8} = N_4 \, \vec{m} \,.
\end{matrix}
\label{eq:xfembasis}
\end{equation}
where $\{N_a\}$ are the underlying shape functions of the continuous element.  One can imagine that the enriched displacement field is represented with a phantom element superimposed over the original.  This is a standard XFEM enrichment, with the slight modification that the weighting coefficients---with support at the nodes of the element---are expressed in the local $\{\vec{n},\vec{m}$\} coordinate system.   This will reduce the number of non-zero entries in the system matrix when coupling with hydrodynamic properties is taken into account. For curving fractures, however, a global coordinate system must be used to avoid ambiguity in defining the local normal and tangent at a node.  In this case, one can simply take $\vec{\xi}_1 = \vec{\eta}_1$, $\vec{\xi}_2 = \vec{\eta}_2$, and so on.  Unlike the EFEM approach, for the XFEM approach inter-element continuity of the enriched displacement field will be enforced below by making the enriched nodal weights global, rather than element-local, degrees of freedom.  

$\;$ 

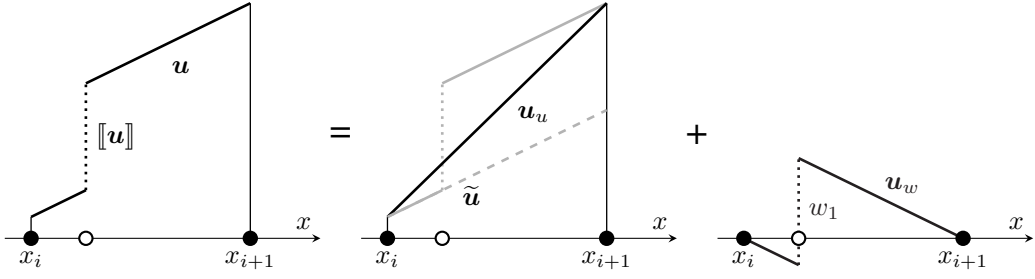
\begin{figure}[t]
	\centering
	\begin{tikzpicture}
	\begin{axis}
	[name = Total,
	width = .35\textwidth,
	height = .35\textwidth,
	xmin=0,   xmax=1.2,
	ymin=-0.13,   ymax=1.1,
	axis y line= none,
	axis x line*=bottom,
	axis lines = middle,
	enlargelimits = true,
	xlabel = {$x$},
	ylabel = {$u$},
	ticks = none
	];
	\addplot [solid, color = black, line width = 1pt] coordinates { (0,0.1) (0.25,0.225) };
	\addplot [solid, color = black, line width = 1pt] coordinates { (0.25,0.725) (1,1.1) };
	\addplot [solid, color = black, line width = .5pt] coordinates { (0.0,0.0) (0.0,0.1) };
	\addplot [solid, color = black, line width = .5pt] coordinates { (1.0,0.0) (1.0,1.1) };
	\draw [dotted, color = black, line width = 1pt] (axis cs:0.25,0.225) -- node[right]{$\llbracket\vec{u}\rrbracket$} (axis cs:0.25,0.725);
	\node[anchor = west, color = black] at (axis cs:0.6, 0.8) {$\vec{u}$};
	\node[anchor = north] at (axis cs:0., -0.02) {$x_i$};
	\node[anchor = north] at (axis cs:1., -0.02) {$x_{i+1}$};
	\draw [thick, fill = black] (axis cs:0., 0.) circle [radius=2.5pt];
	\draw [thick, fill = black] (axis cs:1., 0.) circle [radius=2.5pt];
	\draw [thick, fill = white] (axis cs:0.25, 0.) circle [radius=2.5pt];
	\end{axis}
	\node[anchor = west] at (Total.east)  (equal)  {\Large =};
	\begin{axis}
	[name = Conforming,
	at=(equal.east), 
	anchor=left of west,
	xmin=0,   xmax=1.2,
	ymin=-0.13,   ymax=1.1,
	width = .35\textwidth,
	height = .35\textwidth,
	axis y line = none,
	axis x line*=bottom,
	axis lines = middle,
	enlargelimits = true,
	xlabel = {$x$},
	ylabel = {$u$},
	ticks = none
	];
	\addplot [solid, color = black!30, line width = 1pt] coordinates { (0,0.1) (0.25,0.225) };
	\addplot [solid, color = black!30, line width = 1pt] coordinates { (0.25,0.725) (1,1.1) };
	\addplot [solid, color = black, line width = 0.5pt] coordinates { (0.0,0.0) (0.0,0.1) };
	\addplot [solid, color = black, line width = 0.5pt] coordinates { (1.0,0.0) (1.0,1.1) };
	\draw [dotted, color = black!30, line width = 1pt] (axis cs:0.25,0.225) -- node[right]{} (axis cs:0.25,0.725);	
	\addplot [solid, color = black, line width = 1pt] coordinates { (0,0.1) (1,1.1) };
	\addplot [dashed, color = black!30, line width = 1pt] coordinates { (0, 0.1) (1, 0.6) };
	\node[anchor = west, color = black] at (axis cs:0.55, 0.58) {$\vec{u}_{u}$};
	\node[anchor = west, color = black] at (axis cs:0.3, 0.2) {$\widetilde{\vec{u}}$};
	\node[anchor = north] at (axis cs:0., -0.02) {$x_i$};
	\node[anchor = north] at (axis cs:1., -0.02) {$x_{i+1}$};
	\draw [thick, fill = black] (axis cs:0., 0.) circle [radius=2.5pt];
	\draw [thick, fill = black] (axis cs:1., 0.) circle [radius=2.5pt];
	\draw [thick, fill = white] (axis cs:0.25, 0.) circle [radius=2.5pt];
	\end{axis}
	\node[anchor = west] at (Conforming.east)  (Plus)  {\Large +};
	\begin{axis}
	[name = NonConf,
	at=(Plus.east), 
	anchor=left of west,
	xmin=0,   xmax=1.2,
	ymin=-0.13,   ymax=1.1,
	width = .35\textwidth,
	height = .35\textwidth,
	axis y line = none,
	axis x line*=bottom,
	axis lines = middle,
	enlargelimits = true,
	xlabel = {$x$},
	ticks = none
	];
	\addplot [solid, color = Black, line width = 1pt] coordinates { (0,0) (0.25,-0.125) };
	\addplot [solid, color = Black, line width = 1pt] coordinates { (0.25,0.375  ) (1,0) };
	\draw [dotted, color = Black, line width = 1pt](axis cs:0.25,-0.125) -- node[right]{$w_1$} (axis cs:0.25,0.375);
	\node[anchor = west, color = Black] at (axis cs:0.6, 0.25) {$\vec{u}_{w}$};
	\node[anchor = north] at (axis cs:0., -0.02) {$x_i$};
	\node[anchor = north] at (axis cs:1., -0.02) {$x_{i+1}$};
	\draw [thick, fill = black] (axis cs:0., 0.) circle [radius=2.5pt];
	\draw [thick, fill = black] (axis cs:1., 0.) circle [radius=2.5pt];
	\draw [thick, fill = white] (axis cs:0.25, 0.) circle [radius=2.5pt];
	\end{axis}
	\end{tikzpicture}
	\caption{Decomposition of a one-dimensional displacement field into a continuous field and a discontinuous enrichment.}
	\label{Fig:DisplacementDecomp}
\end{figure}

\begin{figure}[t]
	\newcommand{\panelsize}{0.24}
	\begin{subfigure}[t]{.9\linewidth}
		\begin{subfigure}[t]{\panelsize\linewidth}
			\centering
			\includegraphics[width=\linewidth]{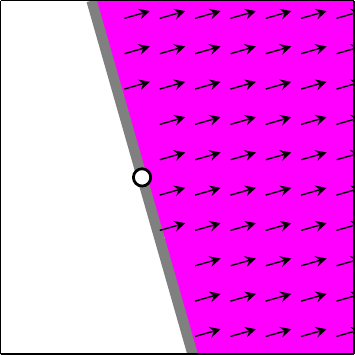}
			\caption{$H\vec{\xi}_1$}
		\end{subfigure} 
		\hfill
		\begin{subfigure}[t]{\panelsize\linewidth}
			\centering
			\includegraphics[width=\linewidth]{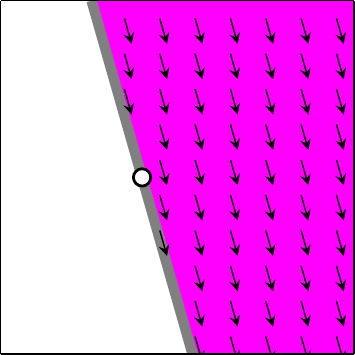}
			\caption{$H\vec{\xi}_2$}
		\end{subfigure}
		\hfill   
		\begin{subfigure}[t]{\panelsize\linewidth}
			\centering
			\includegraphics[width=\linewidth]{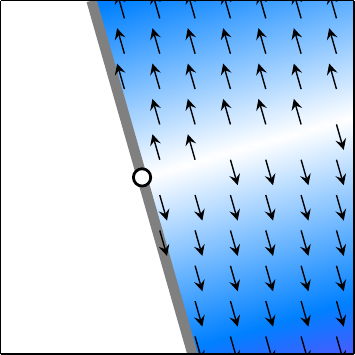}
			\caption{$H\vec{\xi}_3$}
		\end{subfigure}	
		\hfill
		\begin{subfigure}[t]{\panelsize\linewidth}
			\centering 
			\includegraphics[width=\linewidth]{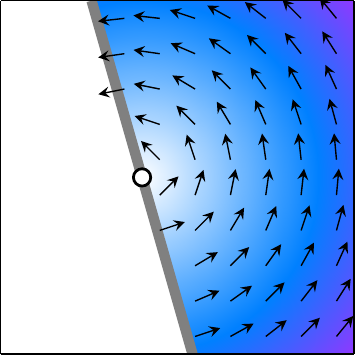}
			\caption{$H\vec{\xi}_4$}
		\end{subfigure}
		
		\vspace{0.2in}
		\begin{subfigure}[t]{\panelsize\linewidth}
			\centering
			\includegraphics[width=\linewidth]{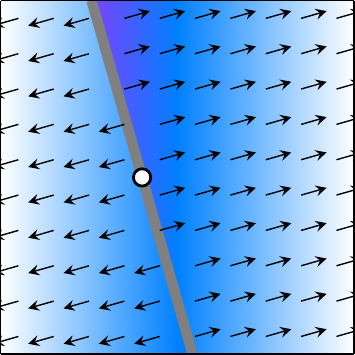}
			\caption{$\vec{\phi}_1$}
		\end{subfigure}\vspace{.7em}
		\hfill
		\begin{subfigure}[t]{\panelsize\linewidth}
			\centering
			\includegraphics[width=\linewidth]{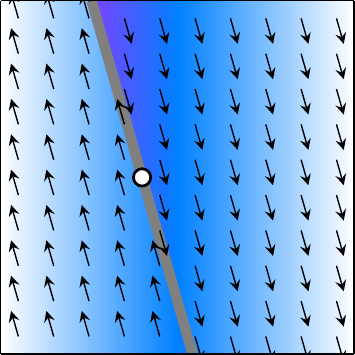}
			\caption{$\vec{\phi}_2$}
		\end{subfigure}\vspace{.7em}
		\hfill
		\begin{subfigure}[t]{\panelsize\linewidth}
			\centering
			\includegraphics[width=\linewidth]{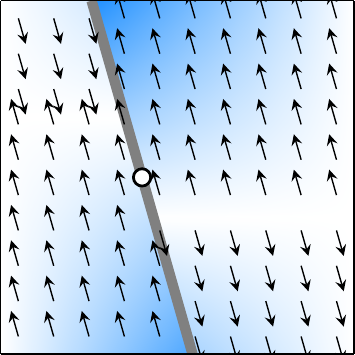}
			\caption{$\vec{\phi}_3$}
		\end{subfigure}\vspace{.7em}
		\hfill
		\begin{subfigure}[t]{\panelsize\linewidth}
			\centering
			\includegraphics[width=\linewidth]{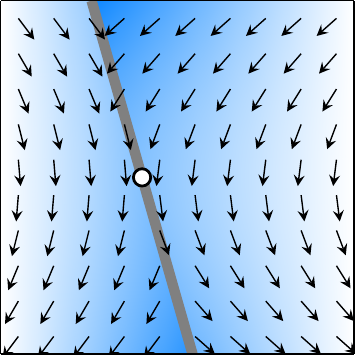}
			\caption{$\vec{\phi}_4$}
		\end{subfigure} 
	\end{subfigure}
	\hfill
	\begin{subfigure}[c]{.06\linewidth}
		\begin{tikzpicture}
		\pgfplotscolorbardrawstandalone[
		colormap/cool,
		point meta min=0,
		point meta max=1,
		parent axis width/.initial=5\linewidth,
		parent axis height/.initial=5\linewidth,
		colormap access=map,
		]
		\end{tikzpicture} 
	\end{subfigure}
	\vspace{-0.2in}
	\caption{Enrichment modes (a-d) and resulting basis functions (e-h) for EFEM(1) in two-dimensions.  The piecewise constant EFEM(0) just uses the first two modes (a-b) and resulting basis functions (e-f). Arrows and colors are representative of direction and magnitude, respectively, of enrichment modes and basis functions.} 
	\label{fig:modes} 
\end{figure}

\begin{figure}[t]
	\newcommand{\panelsize}{0.24}
	\begin{subfigure}[t]{.9\linewidth}
		\begin{subfigure}[t]{\panelsize\linewidth}
			\centering
			\includegraphics[width=\linewidth]{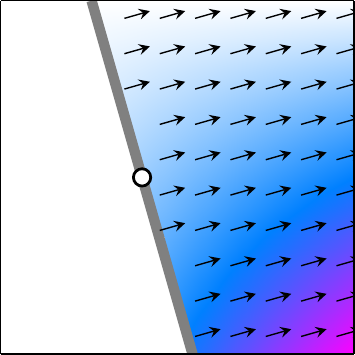}
			\caption{$H\vec{\xi}_1$}
		\end{subfigure} 
		\hfill
		\begin{subfigure}[t]{\panelsize\linewidth}
			\centering
			\includegraphics[width=\linewidth]{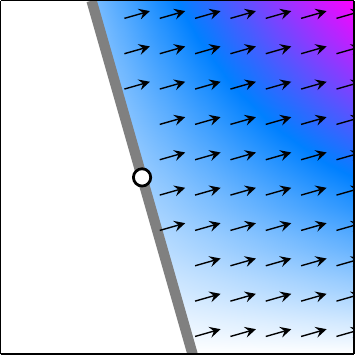}
			\caption{$H\vec{\xi}_2$}
		\end{subfigure}
		\hfill   
		\begin{subfigure}[t]{\panelsize\linewidth}
			\centering
			\includegraphics[width=\linewidth]{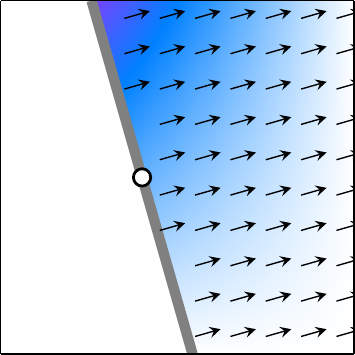}
			\caption{$H\vec{\xi}_3$}
		\end{subfigure}	
		\hfill
		\begin{subfigure}[t]{\panelsize\linewidth}
			\centering 
			\includegraphics[width=\linewidth]{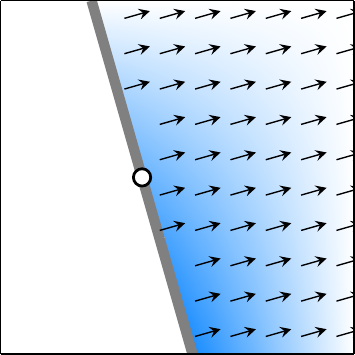}
			\caption{$H\vec{\xi}_4$}
		\end{subfigure}
		
		\vspace{0.2in}
		\begin{subfigure}[t]{\panelsize\linewidth}
			\centering
			\includegraphics[width=\linewidth]{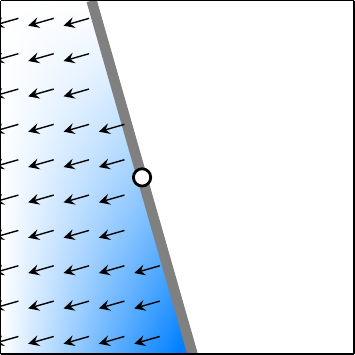}
			\caption{$\vec{\phi}_1$}
		\end{subfigure}\vspace{.7em}
		\hfill
		\begin{subfigure}[t]{\panelsize\linewidth}
			\centering
			\includegraphics[width=\linewidth]{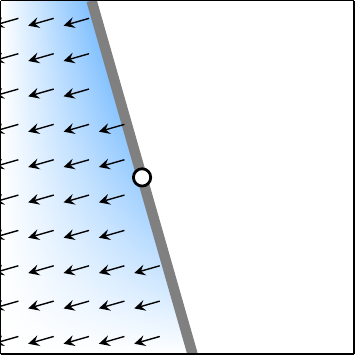}
			\caption{$\vec{\phi}_2$}
		\end{subfigure}\vspace{.7em}
		\hfill
		\begin{subfigure}[t]{\panelsize\linewidth}
			\centering
			\includegraphics[width=\linewidth]{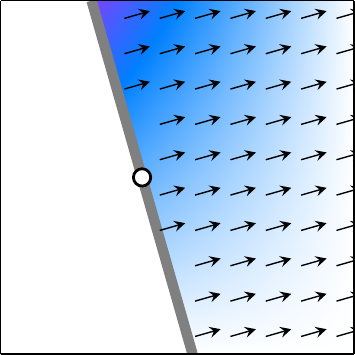}
			\caption{$\vec{\phi}_3$}
		\end{subfigure}\vspace{.7em}
		\hfill
		\begin{subfigure}[t]{\panelsize\linewidth}
			\centering
			\includegraphics[width=\linewidth]{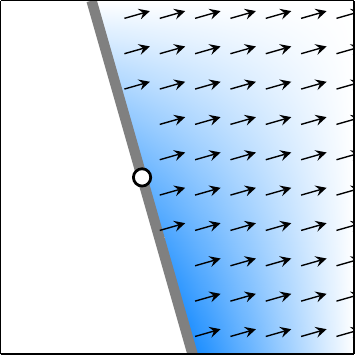}
			\caption{$\vec{\phi}_4$}
		\end{subfigure} 
	\end{subfigure}
	\hfill
	\begin{subfigure}[c]{.06\linewidth}
		\begin{tikzpicture}
		\pgfplotscolorbardrawstandalone[
		colormap/cool,
		point meta min=0,
		point meta max=1,
		parent axis width/.initial=5\linewidth,
		parent axis height/.initial=5\linewidth,
		colormap access=map,
		]
		\end{tikzpicture} 
	\end{subfigure}
	\vspace{-0.2in}
	\caption{Same as Figure \ref{fig:modes} for the four XFEM enrichment modes and basis functions associated with local direction $\vec{n}$ in two-dimensions. Modes and basis associated with the tangential direction have equal magnitude contour plots and direction determined by $\vec{m}$.} 
	\label{fig:xfemmodes} 
\end{figure}

\noindent Regardless of the basis choice, the heaviside-based decomposition equation (\ref{eqn:hdecomp}) is awkward to implement because $\widetilde{\vec{u}} \notin \boldsymbol{U}$.  This field is continuous within elements but discontinuous at element boundaries.  This difficulty may be remedied through a simple manipulation, aimed at expressing the displacement field as sum of a globally continuous function $\vec{u}_u$ and a discontinuous enrichment function $\vec{u}_w$.   Figure \ref{Fig:DisplacementDecomp} illustrates the idea in the simpler one-dimensional case.   Let $\vec{x}_a$ denote the nodal coordinate providing support for shape function $\vec{\eta}_a$.  Also, let $\xi_a$ denote an $x$- or $y$-component of $\vec{\xi}(\vec{x}_a)$---that is, a component of the enriched displacement field evaluated at the node, ordered in the same manner as $\{u_a\}$.  Finally,  we define the function,
\begin{equation}
\vec{f}(\vec{x}) = \sum_{a=1}^8  \xi_a \, H(\vec{x}_a)  \, \vec{\eta}_a (\vec{x}) \,,
\end{equation}
which is a continuous function on $K$.  We proceed by adding and subtracting $\vec{f}$ from the original decomposition,
\begin{align}
\vec{u}_{\left| K \right.} &= (\, \widetilde{\vec{u}}+ \vec{f} \,) + (\,H \, \vec{\xi} - \vec{f} \,) \,, \notag \\
& = \sum_{a=1}^{8} u_a \,\vec{\eta}_a + \sum_{b=1}^{\bar{n}_w} w_b \, \vec{\phi}_b \, ,
\end{align}
where $\bar{n}_w$ is the local number of enrichment modes considered. Here, the new coefficient $u_a = \widetilde{u}_a + H(\vec{x}_a)\, \xi_a $ represents the \emph{total} displacement field evaluated at the node.  This expansion leads to new functions $\{\vec{\phi}_b\}$ with
\begin{equation}
\vec{\phi}_b(\vec{x}) = H(\vec{x}) \, \vec{\xi}_b(\vec{x}) - \sum_{a=1}^{4}   H(\vec{x}_a) \, N_a \vec(\vec{x}) \, \vec{\xi}_b (\vec{x}_a)  \qquad b=1\,, ... \,,\bar{n}_w \,
\end{equation}
It is these functions that form a basis for the space $\boldsymbol{\mathbb{E}}(K)$ on an individual element.  Figures~\ref{fig:modes} and \ref{fig:xfemmodes} present the resulting enrichment functions for EFEM and XFEM in 2D, respectively. By design, these functions are equal to zero at the nodes of the element, but their traces on element edges intersected by a fracture segment are non-zero.  For the EFEM approaches, these traces will not conform between neighbors.  For XFEM, if the same enrichment degree-of-freedom is used for interpolating neighbor elements, edge continuity is automatically enforced.

The global displacement field is now 
\begin{align}
\vec{u} &= \vec{u}_u + \vec{u}_w = \sum_{a=1}^{n_u} u_a \,\vec{\eta}_a + \sum_{b=1}^{n_w} w_b \, \vec{\phi}_b \,,
\label{eq:globaldisp}
\end{align}
with $\vec{u}_u \in \boldsymbol{U}$ and $\vec{u}_w \in \boldsymbol{W}$ as desired.  Here, $n_u$ is the total number of nodal displacement degrees of freedom, and $n_w$ is the total number of element (or nodal) enrichments.  The corresponding strain field is
\begin{equation}
\nabla^s \vec{u} = \sum_{a=1}^{n_u} u_a \, \nabla^s \vec{\eta}_a + \sum_{b=1}^{n_w} w_b \, \nabla^s \vec{\phi}_b \,.
\label{eq:globalstrain}
\end{equation}
Noting that $\nabla H = \delta \vec{n}$, with $\delta$ a Dirac delta function centered on the fracture, we may directly compute the enhanced strain basis functions as 
\begin{equation}
\nabla^s \vec{\phi}_b = H \nabla^s \vec{\xi}_b
- \nabla^s \vec{f} + \delta (\vec{\xi}_b \otimes \vec{n})^s \,.
\label{eq:strainbasis}
\end{equation}
The first two terms in the strain are regular contributions, while the Dirac delta term is singular and is only present on the fracture segment $k$.  Note that the effective stress in the matrix is a direct function of the regular strain, which includes a continuous and enhanced contribution.

We may now proceed to the variational form of the mechanical problem.  The EFEM schemes are of non-Galerkin type, in which the space of test strains $\nabla^s \widehat{\vec{v}}$ will differ from the space of trial strains $\nabla^s \vec{u}$.  We begin, however, with the Galerkin approach---as typically used for XFEM---and consider trial functions $\vec{v} = \vec{\eta}+\vec{\phi}$,  with $\vec{\eta} \in \boldsymbol{U}$ and $\vec{\phi} \in \boldsymbol{W}$.
The discrete variational form of the momentum balance at time $t_\tau$ is to find $\{\vec{u}_u,\vec{u}_w\}_\tau$  such that
\begin{equation}
\int_\mathcal{B} \nabla^s \vec{v} : \vec{\sigma}_\tau \,dA - \int_\mathcal{B} \vec{v} \cdot \rho_\tau \vec{g} \,dA - 
\int_{\partial \mathcal{B}_t} \vec{v} \cdot \overline{\vec{t}}_\tau \,dL 
= \vec{0}
\qquad
\forall \{\vec{\eta},\vec{\phi}\} \,.
\end{equation}
Note that the total stress and total traction can be expanded into mechanical and fluid contributions, but for the moment it is convenient to work in the more compact notation of total quantities.  Using equations (\ref{eq:globaldisp}--\ref{eq:strainbasis}), the test functions may be expanded, and Dirac delta contributions converted from area to line integrals.  The independence of the variations then leads to the following weak form: Find $\{ \vec{u}_u,\vec{u}_w \}_\tau \in  \boldsymbol{U} \times  \boldsymbol{W}$ such that
\begin{align}
\mathcal{R}^u &= \int_\mathcal{B} \nabla^s \vec{\eta} : \vec{\sigma}_\tau \,dA - \int_\mathcal{B} \vec{\eta} \cdot \rho_\tau \vec{g} \,dA - 
\int_{\partial \mathcal{B}_t} \vec{\eta} \cdot \overline{\vec{t}}_\tau \,dL = \vec{0} \,,   \label{eq:galerkin_mom}\\
\mathcal{R}^w &= \int_\mathcal{M} \nabla^s \vec{\phi} : \vec{\sigma}_\tau \,dA - \int_\mathcal{M} \vec{\phi} \cdot \rho_\tau \vec{g} \,dA - 
\int_\mathcal{F} \vec{\xi} \cdot {\vec{t}}_\tau \,dL   = \vec{0} \,,
\label{eq:galerkin_traction}
\end{align}
for all $\{ \vec{\eta},\vec{\phi} \} \in  \boldsymbol{U} \times  \boldsymbol{W}$.
Here, $\vec{t}_\tau$ is the local traction on the fracture surface, evaluated as a function of the local displacement jump. The first balance leads to $n_d$ discrete residual equations enforcing a global momentum balance.  These equations are coupled through the nodal support of the basis functions.  The second residual leads to $n_w$ discrete equations that enforce traction balance across fracture segments.  In the XFEM approach, these additional enrichments also have nodal support and therefore couple neighbor elements if they are both enriched.  In EFEM, the enrichments are local bubbles, and equations (\ref{eq:galerkin_traction}) can be reduced to $n_w$ element-wise equations.  This has implications for the system matrix sparsity and resulting linear solver strategy.

Unfortunately, due to the reduced kinematics used in the EFEM schemes, the traction continuity within cut elements may be poorly approximated in certain configurations.  See \cite{Jirasek2000,Wu2011} for an extensive discussion on this topic.  A non-Galerkin formulation is widely preferred for improving these deficiencies.  Applying the divergence theorem to (\ref{eq:galerkin_traction}) implies,
\begin{equation}
\int_\mathcal{M}  \vec{\phi}_b \cdot \left( \nabla \cdot \vec{\sigma}_\tau +  \rho_\tau \vec{g} \right) \,dA +
\int_\mathcal{F} \vec{\xi}_b \cdot \left( \vec{\sigma}_\tau \cdot \vec{n} - {\vec{t}}_\tau \right) \,dL 
= \vec{0}\,,  \qquad b=1,2,... \, ,n_w \,.
\end{equation} 
If one assumes the balance of linear momentum is satisfied in strong form as the method converges, the first term on the left-hand side vanishes.  The remaining equation directly expresses a weak enforcement of traction continuity, as desired.  Unfortunately, stresses are typically evaluated at the quadrature points in the bulk of the element, not on the fracture surface.  To remedy this, let us introduce a projection operator that maps stresses in the bulk to tractions on the surface, 
\begin{align}
\int_\mathcal{M} ( \vec{\beta}_b \otimes \vec{n})^s : \vec{\sigma}_\tau \,dA - 
\int_\mathcal{F} \vec{\xi}_b \cdot {\vec{t}}_\tau \,dL 
& = \vec{0} \,,
\qquad b=1,2,... \, ,n_w \,,
\end{align}
with the projection operator $( \vec{\beta}_b \otimes \vec{n})^s$ playing the role of a new test strain.  This is equivalent to writing
\begin{align}
\int_\mathcal{M} \vec{\beta}_b \cdot \vec{t}_\tau \,dA - 
\int_\mathcal{F} \vec{\xi}_b \cdot {\vec{t}}_\tau \,dL 
& = \vec{0} \,,
\qquad b=1,2,... \, ,n_w \,,
\label{eq:modbalance}
\end{align}
highlighting the fact that this form directly imposes an equilibrium between tractions evaluated on the fracture segment and traction values $\vec{t}_n = \vec{\sigma}_\tau \cdot \vec{n}$ evaluated in the element volume.  

The key question is how to choose the test vectors $\vec{\beta}_b$. One strategy, proposed in \cite{Linder2007} for the EFEM(1) basis, is to determine the vectors $\{\vec{\beta}_b\}$  such that equation~(\ref{eq:modbalance}) will be \emph{exactly} satisfied for any traction field that is a piecewise polynomial up to a certain order.  For the EFEM(0) interpolation, for example, let us assume this balance should be satisfied for a piecewise constant traction field in each cut element,
\begin{align}
\vec{t}_{\left| K \right.} & = t_1 \vec{n} + t_2 \vec{m} \,.
\label{eq:consttraction}
\end{align}
The test basis is also chosen as piecewise constant,
\begin{equation}
\vec{\beta}_{1 \left| K \right.} = \beta_1 \vec{n} \qquad \vec{\beta}_{2\left| K \right.} = \beta_2 \vec{m}\,,
\label{eq:consttestbasis}
\end{equation}
with two unknown coefficients $\beta_1,\beta_2$. Restricting equation~(\ref{eq:modbalance})  to an individual element, inserting expressions (\ref{eq:constbasis}), (\ref{eq:consttraction}), and (\ref{eq:consttestbasis}), and noting that the traction coefficients are arbitrary implies $\beta_1 = \beta_2 = |k|/|K|$.  That is, the scaling weight is the ratio of the fracture length to the element area. Because the normal and tangential tractions are orthogonal, the weights in both directions are equal. 
In the case of the EFEM(1) discretization, we assume the traction is a piecewise linear polynomial on a cut element,
\begin{align}
\vec{t}_{\left| K \right.} & = \left(t_1 + t_2 y_n + t_3 y_m \right) \vec{n} + \left(t_4 +  t_5y_n + t_6 y_m \right) \vec{m}
\label{eq:lintraction}
\end{align}
The four components of the test basis are chosen as linear polynomials,
\begin{align}
\vec{\beta}_{1\left| K \right.} &= \left(\beta_1 + \beta_2 y_n + \beta_3 y_m \right)\vec{n} \qquad
\vec{\beta}_{2\left| K \right.} = \left(\beta_4 + \beta_5 y_n + \beta_6 y_m \right)\vec{n} \notag \\
\vec{\beta}_{3\left| K \right.} &= \left(\beta_1 + \beta_2 y_n + \beta_3 y_m \right)\vec{m} \qquad
\vec{\beta}_{4\left| K \right.} = \left(\beta_4 + \beta_5 y_n + \beta_6 y_m \right)\vec{m} 
\label{eq:betavecs}
\end{align}
with six weighting coefficients $\{\beta_i\}$ to be determined.  Restricting equation \eqref{eq:modbalance} to an individual element, inserting relationships \eqref{eq:linderbasis}, \eqref{eq:lintraction}, and \eqref{eq:betavecs}, and noting that the balance must be satisfied for arbitrary $\{t_i\}$ implies
\begin{equation}
\begin{bmatrix}
\beta_1 & \beta_4 \\
\beta_2 & \beta_5 \\
\beta_3 & \beta_6 
\end{bmatrix}
= \Mat{A}^{-1}_K \Mat{L}_k
\end{equation}
with matrices
\begin{equation}
\Mat{A}_K = \int_K 
\begin{bmatrix}
1 & y_n & y_m 
\end{bmatrix}^T
\begin{bmatrix}
1 & y_n & y_m 
\end{bmatrix} 
dA
\qquad \text{and} \qquad
\Mat{L}_k = \int_k 
\begin{bmatrix}
1 & 0 & y_m 
\end{bmatrix}^T
\begin{bmatrix}
1 &  y_m 
\end{bmatrix} 
dL\,.
\end{equation}
The $3\times3$ matrix $\Mat{A}_K$ is a mass matrix for the $\mathbb{P}_1(K)$ basis $\{1,y_n,y_m\}$.  Note that $y_n = 0$ on $k$ because the coordinate system origin is set on the fracture surface.  The test vector weights only depend on the cut element geometry and may be computed in a pre-processing step.

In summary, in the EFEM methods the global momentum balance equation~(\ref{eq:galerkin_mom}) remains the same, but the discrete traction balance residual equations at time $t_\tau$ are replaced as: 
%
\begin{align}
\widetilde{\mathcal{R}}^{w}_{b} =&\;
\int_\mathcal{M} ( \vec{\beta}_b \otimes \vec{n})^s : \vec{\sigma}_\tau \,dA - 
\int_\mathcal{F} \vec{\xi}_b \cdot {\vec{t}}_\tau \,dL 
= \vec{0}  \,,
\qquad b=1,2,... \, ,n_w \,,
\label{eq:mech_residuals_b} 
\end{align}
Note that the symmetry of the Galerkin form in XFEM is lost in the EFEM methods.  We also emphasize that the total stress,  total traction, and density terms contain both solid and fluid contributions, and therefore these residual equations are tightly coupled to all of the unknown fields.


\begin{remark}
	This section began with the assumption that elements are cut by at most one fracture segment.  In geologic media, however, this assumption is highly restrictive.  Intersecting fracture networks are pervasive in the subsurface.  In the numerical examples below, when an element is cut by multiple fracture segments---either in a sub-parallel fashion or fully intersecting---we adopt a simple superposition treatment of the problem.  The resulting strain field is written as the additive sum of the separate fracture contributions. For example, for an element cut by two fractures, the total strain field reads
	\begin{equation}
	\label{eq:intersectionStrain}
	\nabla^s \vec{u} = \sum_{a=1}^{n_u} u_a \, \nabla^s \vec{\eta}_a + \sum_{b=1}^{n_w} w^{(1)}_b \, \nabla^s \vec{\phi}^{(1)}_b +\sum_{b=1}^{n_w} w^{(2)}_b \, \nabla^s \vec{\phi}^{(2)}_b \, ,
	\end{equation}
	where the superscript $(i)$ indicates the quantities relative to fracture $i$. Note, that a traction balance, as expressed in equation \eqref{eq:modbalance}, has to be considered for each fracture segment and that a local coupling between the unknowns relative to each fracture is introduced.
	This assumption is clearly a significant simplification of the kinematics, but it retains the simplicity of the underlying method.  Further, in practical geologic applications one often observes that significant slip only occurs on a subset of well-oriented fractures, while poorly oriented fractures have weaker interactions.  Interesting work on more complicated enrichments to directly treat fracture intersections can be found in \cite{Linder2009}. 
\end{remark}

\begin{remark}
	In the XFEM scheme, the test strain $\nabla^s \vec{\phi_b}$ is discontinuous, and standard element quadrature will be inaccurate when evaluating equation (\ref{eq:galerkin_traction}).  Here, we sub-triangulate cut elements for quadrature purposes.  In the EFEM schemes, standard Gaussian quadrature may be applied.
\end{remark}

\begin{remark}
	In the EFEM(1) scheme, the enriched strain can have linearly dependent columns when a fracture cuts only one node of an element.  This can lead to singularity of the system matrix.  In \cite{Linder2007}, a stabilization is recommended in which the higher-order modes are penalized so that the scheme more closely mimics the intrinsically stable EFEM(0) scheme on such elements.   This is one downside of the EFEM(1) scheme, as single-node-cut elements are difficult to avoid for arbitrarily oriented fractures.
\end{remark}

\begin{remark}
	In XFEM, the conditioning of the system matrix is sensitive to the ratio of element areas bisected by the fracture.  In particular, fracture surfaces that pass close to element nodes may cause ill-conditioned or even singular matrices.  Again, such a configuration is difficult to avoid for arbitrarily oriented fractures.  In \cite{Babuvska2012, Wu2015}, stabilized versions of XFEM are proposed to address this issue.  In particular, the method in \cite{Wu2015} appears to provide stable matrix conditioning while maintaining the same approximation accuracy as the basic method.
\end{remark}

\begin{remark}
	For simplicity, we have not included any special treatment of fracture tips.  It is well known, however, that numerical error due to singular strain fields at the tips can dominate convergence behavior.  Tip enrichments are commonplace in XFEM implementations, but comparable equivalents for EFEM methods are not widely used.  It has been observed, however, that the singular nature of the non-conforming jump at the tip in the EFEM scheme can partially compensate for the lack of a dedicated tip enrichment \cite{Borja2008}.  Nevertheless, the development of effective tip treatments for EFEM methods for use in hydraulic fracturing applications is the subject of ongoing work. 
\end{remark}

\begin{remark}
	In XFEM, one can ensure that the fracture aperture smoothly closes in a tip element by only enriching two of the four nodes.  In the numerical results below, however, we found better accuracy by allowing for a discontinuous jump at the tip by enriching all four nodes.  This approach mimics the EFEM treatment, allowing for a non-conforming tip jump.
\end{remark}

\begin{remark}
	In this work, fractures are assumed to be stationary. In some practical scenarios, however, fractures can grow due to hydraulic stimulation. The use of an embedded method then has the great advantage of allowing fractures to propagate in non-grid-aligned orientations. Neverthess, non-trivial challenges remain in using EFEM \cite{Armero2009,ArmeroFracProp2012,Linder2013} and XFEM \cite{moes1999finite,sukumar2000extended,Khoei2014} for fracture growth, particularly in the evaluation of tip propagation criteria.
\end{remark}

\begin{remark}
	While the small strain assumption is adopted in this work, EFEM has also been used in the literature to deal with discontinuities in presence of finite deformations \cite{ArmerofiniteDef2008}.
\end{remark}

\subsection{Flow Discretization}

We adopt a finite volume scheme for discretizing the flow physics, because of its element-wise mass conservation properties.  The enrichment of the pressure and saturation fields to address fractured elements is therefore straightforward (Figure~\ref{fig:discretization}).  In the matrix, there is a piecewise constant pressure field $p^m(\vec{x})$ and saturation field $s^m(\vec{x})$, each with one degree of freedom per cell $K$.  Similarly, in the fracture network we solve for a piecewise constant pressure field $p^f(\vec{x})$ and saturation field $s^f(\vec{x})$, each with one degree of freedom per fracture segment $k$. In keeping with the finite element formalism used above for the mechanical portion of the problem, we can define two discrete spaces \cite{Barth2018,Eymard2000},
\begin{align}
\mathcal{P}^m & := \left\lbrace \psi^m \left|\right. \psi^m \in L^2(\mathcal{B}), \psi^m_{\left| K \right.} \in \mathbb{P}_0(K) \; \forall K \in \mathcal{T}_\mathcal{B}  \right\rbrace,	\\
\mathcal{P}^f & := \left\lbrace \psi^f \left|\right. \psi^f \in L^2(\mathcal{F}), \psi^f_{\left| k \right.} \in \mathbb{P}_0(k) \; \forall k \in \mathcal{T}_\mathcal{F}  \right\rbrace,	
\end{align}
which define the space of piecewise constant fields on the given triangulation (matrix or fracture).   Letting $\{ \psi^{\alpha}_a \}$ denote a set of basis functions for $\mathcal{P}^\alpha$, $\alpha=\{m,f\}$, the pressure and saturation fields are interpolated as
\begin{equation}
p^\alpha (\vec{x}) =  \sum_{a=1}^{n_\alpha} p_a^\alpha \psi^\alpha_a (\vec{x}) 
\qquad \text{and} \qquad 
s^\alpha (\vec{x}) =  \sum_{a=1}^{n_\alpha} s_a^\alpha \psi^\alpha_a (\vec{x}) \,.
\end{equation}
Note that the basis functions are simply characteristic functions, equal to one on a given volume and zero elsewhere.  The discrete governing equations are then functions of the fields $\{p^m,s^m\} \in [\mathcal{P}^m]^2$ and $\{p^f,s^f\} \in [\mathcal{P}^f]^2$. 

In the flow formulation, it is also necessary to define discrete fluxes between neighboring volumes and from external boundaries.   These fluxes come in four forms: boundary-to-matrix, matrix-to-matrix, fracture-to-fracture, and matrix-to-fracture.  It is easiest to manage these fluxes using connectivity lists.  Let $\mathcal{C}^{mm}$ denote the set of unique pairs $(J,K)$ of mesh cells connected through a common face.  Similarly, let $\mathcal{C}^{ff}$ denote the set of neighboring fracture segments $(j,k)$, $\mathcal{C}^{mf}$ the set of fracture segments embedded in matrix cells $(K,k)$, and $\mathcal{C}^{em}$ a list of external connections $(\bar{j},K)$ allowing for non-zero boundary fluxes from boundary edge $\bar{j}$.

Using a standard two-point flux approximation (TPFA), the matrix-to-matrix ($mm$), fracture-to-fracture ($ff$), or matrix-to-fracture ($mf$) discrete fluxes of phase $\pi$ are computed as
\begin{subequations}\label{eqn:tpfa}
	\begin{align}
	F^{mm}_{\pi,\,JK} &=  - \frac{\rho_\pi^{\text{upw}} k_{r\pi}}{\mu_\pi^{\text{upw}}} T_{JK}[(p^m_K + \varrho_\pi g z^m_K ) - ( p^m_J - \varrho_\pi g z^m_J) ], & (J,K) \in \mathcal{C}^{mm}, \\
	F^{ff}_{\pi,\,jk} &=  - \frac{\rho_\pi^{\text{upw}} k_{r\pi}}{\mu_\pi^{\text{upw}}} T_{jk}[(p^f_k + \varrho_\pi g z^f_k ) - ( p^f_j - \varrho_\pi g z^f_j) ], & (j,k) \in \mathcal{C}^{mm}, \\
	F^{mf}_{\pi,\,Kk} &=  - \frac{\rho_\pi^{\text{upw}} k_{r\pi}}{\mu_\pi^{\text{upw}}} T_{Kk}[(p^f_k + \varrho_\pi g z^f_k ) - ( p^m_K - \varrho_\pi g z^m_K) ], & (K,k) \in \mathcal{C}^{mf}.
	\end{align}
\end{subequations}
The superscript $(\cdot)^\text{upw}$ denotes an upwinded quantity, whereas $\varrho_\pi$ denotes the phase mass density averaged at the interface between control volumes.  For $mm$ fluxes, the transmissibility coefficient between control volumes $J$ and $K$ connected by face $c$ is computed as the harmonic average
\begin{equation}
T_{JK} = \frac{T_{Jc} T_{Kc}}{T_{Jc} + T_{Kc}} \,,
\end{equation}
with half-transmissibility coefficients given by
\begin{equation}
T_{Jc} = |c| \frac{(\vec{x}_c - \vec{x}_J) \cdot \vec{\kappa} \cdot \vec{n}_{J,c}}{\| \vec{x}_c - \vec{x}_J \|^2} 
\qquad \text{and} \qquad
T_{Kc} = |c| \frac{(\vec{x}_c - \vec{x}_K) \cdot \vec{\kappa} \cdot \vec{n}_{K,c}}{\| \vec{x}_c - \vec{x}_K \|^2} \,.
\end{equation}
Here, $|c|$ is the area---namely, the $(d-1)$-measure, with $d$ the spatial dimension of the problem---of the connecting face, $\vec{x}_c$ is a suitably chosen collocation point on the face, and $\vec{x}_{J}$ and $\vec{n}_{J,c}$ (respectively $\vec{x}_{K}$ and $\vec{n}_{K,c}$) are the centroid and outer unit normal at the connecting face for control volume $J$ (respectively $K$).  The same form is used for the segment-to-segment (fracture) transmissibility $T_{jk}$ through a suitable interpretation of the required quantities \cite{MousaFADM2018}.  For example, the area of a connecting fracture face is the aperture times the $(d-2)$-measure of the intersection between boundaries of control volumes $j$ and $k$. The key difference with the matrix flux is that the aperture and absolute permeability of fracture segments is strongly dependent on the mechanical deformation, while for matrix cells the connecting area and absolute permeability is assumed constant. For the case of mass interchange between a fracture element $k$ embedded in cell $K$, the transmissibility coefficient is computed as 
\begin{equation}
T_{Kk} = \frac{|\,k \,|}{\langle d_{Kk} \rangle} \frac{\kappa_K \kappa_k}{\kappa_K + \kappa_k} \,.
\end{equation}
Here, $|\,k \,|$ is the $(d-1)$-measure of fracture that cuts through matrix element $K$, $\kappa_K = (\vec{n} \otimes \vec{n}) : \vec{\kappa}_K$ is a scalar measure of the matrix permeability, and $\left\langle d_{Kk}\right\rangle$ is an average connection distance defined as
%
%
\begin{equation}
\left\langle d_{Kk}\right\rangle = \frac{1}{|\,K \,|} \int_K | (\vec{x}-\vec{x}_k) \cdot \vec{n} |  \, dA \,,
\end{equation}
where $\vec{x}_k$ is the centroid of the fracture segment. 
For external fluxes $F^{em}_{\pi,\,\bar{j}K}$, Neumann fluxes may be directly prescribed.  On a Dirichlet boundary, the flux formula (\ref{eqn:tpfa}) may be used with the transmissibility replaced by the half-transmissibility, and one of the pressures interpreted as a known external pressure $\overline{p}_b$.

The  weak form of the mass balance equations at time $t_\tau$ may then be stated as: Find $\{ p^m, s^m,p^f,s^f \}_\tau \in  [\mathcal{P}^m]^2 \times  [\mathcal{P}^f]^2$ such that, for each phase $\pi= \{w,nw\}$,
%
\begin{align}
\mathcal{R}^m_\pi =&\;
\int_{\mathcal{B}} \psi_\pi^m \left( \frac{\Delta m^m_{\pi} }{\Delta t} - q^m_{\pi,\tau} \right) \, \mathrm{d}A -
\sum_{C^{mm}} (\psi^m_{\pi,K}-\psi^m_{\pi,J}) F_{\pi,\tau,JK}^{mm} 
- \sum_{C^{mf}} (\psi^f_{\pi,k}-\psi^m_{\pi,K}) F_{\pi,\tau,Kk}^{mf} 
+ \sum_{C^{em}} \psi^f_{\pi,K} F_{\pi,\tau,\bar{j}K}^{em}    = 
0   \label{eq:flow_residuals_a} \\
\mathcal{R}^f_\pi =&\;
\int_{\mathcal{F}} \psi_\pi^f \left( \frac{\Delta m^f_{\pi} }{\Delta t} - q^f_{\pi,\tau} \right) \, \mathrm{d}A -
\sum_{C^{ff}} (\psi^f_{\pi,k}-\psi^f_{\pi,j}) F_{\pi,\tau,jk}^{ff} 
+ \sum_{C^{mf}} (\psi^m_{\pi,k}-\psi^f_{\pi,K}) F_{\pi,\tau,Kk}^{mf}    = 
0   \label{eq:flow_residuals_b} 
\end{align}
for all $\{ \psi^m_w, \psi^m_{nw},\psi^f_w,\psi^f_{nw} \} \in  [\mathcal{P}^m]^2 \times  [\mathcal{P}^f]^2$. We emphasize that the above residuals represent four discrete equations, two for each phase in the matrix and fractures.  Because the basis functions are simple characteristic functions, these equations may be readily assembled using an element based accumulation loop followed by a connection based flux loop, with appropriate indexing and signs. 

\begin{remark}
	This formulation only provides accurate results for conductive fractures. If fractures have permeabilities lower than the rock matrix, and are thus flow barriers, the p-EDFM formulation provides an alternative \cite{Tene2017,Jiang2017}.
\end{remark}

\subsection{Linearization and solution strategy}
Equations~(\ref{eq:galerkin_mom}, \ref{eq:mech_residuals_b}, \ref{eq:flow_residuals_a}, \ref{eq:flow_residuals_b}) lead to a set of algebraic residual equations describing the coupled behavior of the system. n this work, this system of equations is solved following a fully-coupled strategy. Let $\blkVec{x}_\tau$ denote a vector gathering all of the nodal and cell-based unknowns at time $t_\tau$.  The residual system may be compactly written as
\begin{equation}
\blkVec{r} (\blkVec{x}_\tau,\blkVec{x}_{\tau-1}) = \blkVec{0} \,.
\end{equation}
The nonlinear equations are solved using Newton's method, with a backtracking algorithm to improve convergence robustness.  Given an estimate of the new solution $\blkVec{x}^{i-1}_\tau$ at time $\tau$ and iteration $(i-1)$, an improved estimate $\blkVec{x}^{i}_\tau$ is determined by
\begin{align}
\text{solving} \quad & \blkMat{J}_\tau^{i-1} \Delta \blkVec{x} = - \blkVec{r}^{i-1}_\tau \,, \\
\text{updating} \quad & \blkVec{x}^i_\tau = \blkVec{x}^{i-1}_\tau + \gamma \Delta \blkVec{x} \,.
\end{align}
Here, $\blkMat{J} = \partial \blkVec{r} / \partial \blkVec{x}$ is the Jacobian of the nonlinear system, and $\gamma \in (0,1]$ is a backtracking parameter that limits the step length in direction $\Delta \blkVec{x}$ to ensure a residual reduction.  The iterations are terminated when the residual norm drops below a desired convergence tolerance, $\| \, \blkVec{r} \, \| < \text{tol}$.  

Given the coupling within the underlying PDEs, the Jacobian system has the block structure,
\begin{equation}
\underbrace{
	\left[\begin{array}{cccccc}
	\square & \square & \square && \square & \\
	\square & \square & \square & \square \\
	\square & \square & \square & \square & \square & \square\\
	& \square & \square & \square & \square & \square\\ 
	\square & \square & \square & \square & \square & \square\\
	& \square & \square & \square & \square & \square\\ 
	\end{array} \right]
}_{\blkMat{J}}
\underbrace{
	\begin{bmatrix}
	\Delta \Vec{u}\\
	\Delta \Vec{w}\\ 
	\Delta \Vec{p}^{m}\\ 
	\Delta \Vec{p}^{f}\\ 
	\Delta \Vec{s}^{m}\\  
	\Delta \Vec{s}^{f}
	\end{bmatrix}
}_{\Delta \blkVec{x}}
= -\underbrace{
	\begin{bmatrix}
	\Vec{r}\sub{u}\\
	\Vec{r}\sub{w}\\ 
	\Vec{r}^m\sub{o}\\ 
	\Vec{r}^f\sub{o}\\ 
	\Vec{r}^m\sub{w}\\ 
	\Vec{r}^f\sub{w}
	\end{bmatrix}
}_{\blkVec{r}} \,,
\label{eq:linsys}
\end{equation}
where $\square$ indicates a sparse matrix defining a non-zero coupling between the respective fields.  For brevity, we have omitted the specific forms of the linearized operators.  They may be directly---if somewhat tediously---derived from the residual equations in a standard way.

\begin{remark}
	In this work, the linear system~(\ref{eq:linsys}) is handled using a direct solver.  This is sufficient for two-dimensional or small three-dimensional problems, but in general it will not provide a scalable approach.  The design of effective preconditioned iterative methods for this system is the subject of ongoing work, building on techniques presented in \cite{White2019,Bui2020}.
\end{remark}

\section{Numerical experiments}

\label{Sec:NumResults}
We now present several numerical experiments to illustrate the relative performance of the schemes.  All rock and fluid properties employed in the simulations are summarized in Table \ref{Tab:TestCasesProperties}.
\begin{testcase}
	A single fracture, embedded in an infinite medium, is pressurized and opens. This problem has an analytical solution\cite{Sneddon1946}.  It is employed to study sensitivity to enrichment strategy, grid resolution, and grid orientation.
\end{testcase}
\begin{testcase}
	We reproduce a test case from \cite{Borja2008} involving a single fracture subject to compression. As such, the fracture slips and provides an opportunity to validate the schemes under shearing. 
\end{testcase}
\begin{testcase}
	We reproduce a test case from \cite{Rethore2007,Khoei2014} involving coupled single-phase flow and geomechanics.  Our results are compared to independent results obtained using an XFEM-based scheme in the original references.
\end{testcase}
\begin{testcase}
	A vertical section of a fractured heterogeneous reservoir is considered. The domain is saturated with a viscous phase which is extracted by a production well. Production is stimulated by injection of a slightly less viscous fluid phase. It includes much of the complexity encountered in multiphase field applications. As such, it is a proof-of-concept for the applicability of the method to realistic engineering scenarios.
\end{testcase}

\begin{table}[htb]
	\caption{Rock and fluid properties. }
	\centering
	\begin{tabular}{lllrrrr}
		\toprule
		&&& \bf Test Case 1 & \bf Test Case 2 & \bf Test Case 3 & \bf Test Case 4 \\
		\midrule
		$K_{dr} $ & Drained skeleton modulus &[GPa] & 15 & $8.3\times 10^{-5}$ & 15 & 11.3\\
		$\nu$ &Poisson's ratio &[-] & 0.4 & 0.3 & 0.4 & 0.25\\
		$\alpha$ & Friction coefficient & [-] & - & 0.1  & 0.6  & 0.6 \\
		$b$ & Biot's coefficient & [-] & 1 & 1 & 1 & 1\\
		$\phi_0$ &Reference porosity &[-] & - & - & 0.3 & 0.3\\
		$\vec{\kappa}$ & Permeability &[$\text{m}^2$] & - & - & $10^{-12}$ & Fig. \ref{Fig:Case4_GeometryPerm} \\
		$K_{w}$ &Water bulk modulus &[GPa] & - & - & $10^{27}$& $10^9$ \\
		$\mu_w$ &Water viscosity& [$\text{Pa}\, \text{s}$] & - & - &  $10^{-3}$ & $10^{-3}$ \\
		$\rho_{w}$ &Water density& [${\text{kg}}/{\text{m}^3}$] & - & - & $1000$ & $1000$\\
		$K_{o}$& Oil bulk modulus& [GPa] & - & - & - & $10^9$ \\
		$\mu_o$ &Oil viscosity &[$\text{Pa}\, \text{s}$] & - & - & - & $1.5 \times10^{-3}$ \\
		$\rho_{o}$& Oil density& [${\text{kg}} / {\text{m}^3}$] & - & - & - & $850$\\
		\bottomrule
	\end{tabular}
	\label{Tab:TestCasesProperties}
\end{table}     

\subsection{Test Case 1 - Opening of a single fracture in an infinite medium}
To begin, a single fracture under constant fluid pressure in an infinite 2D medium is considered. An analytical solution to this problem is available \cite{Sneddon1946}.  In this test case, no fluid flow occurs, so only the mechanical equations are solved.  For the numerical solution, a $100 \, \text{m} \times 100 \, \text{m} $ square domain is considered and a fracture is placed at several rotation angles in the center. A $15 \times 15$ cartesian grid is used as a base mesh, and uniform  $3\times3$ refinements are performed 3 times to study error convergence behavior.  The initial fracture length at a given angle is chosen so that the fracture terminates at an element boundary. The fracture is pressurized at a constant pressure $p^f = 1.0$ MPa. To avoid boundary effects, analytical displacement boundary conditions are imposed on the external boundaries.  The boundary solution is numerically computed using a displacement discontinuity method for an infinite body as described in \cite{Crouch1976}. 

\begin{figure}[htbp]
	\newcommand{\panelsize}{0.275}

	\begin{subfigure}[t]{\panelsize\textwidth}  
		\centering
		\includegraphics[width=\linewidth]{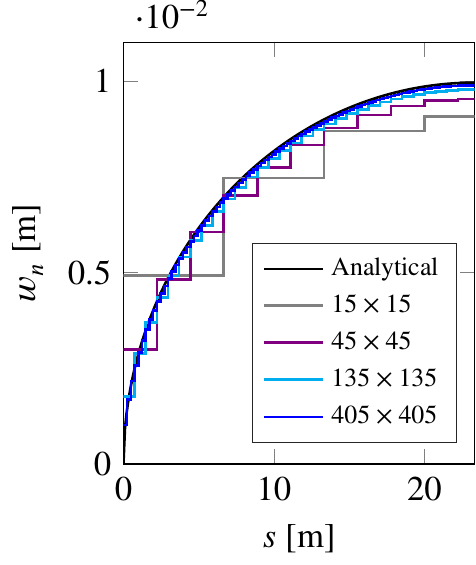}
		\caption{EFEM(0)}
	\end{subfigure}	
	\hfill	
	\begin{subfigure}[t]{\panelsize\textwidth}  
		\centering
		\includegraphics[width=\linewidth]{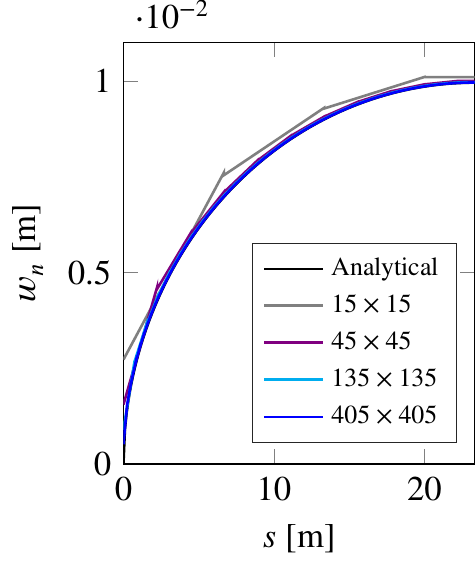}
		\caption{EFEM(1)}
	\end{subfigure}	
	\hfill	
	\begin{subfigure}[t]{\panelsize\textwidth}  
		\centering
		\includegraphics[width=\linewidth]{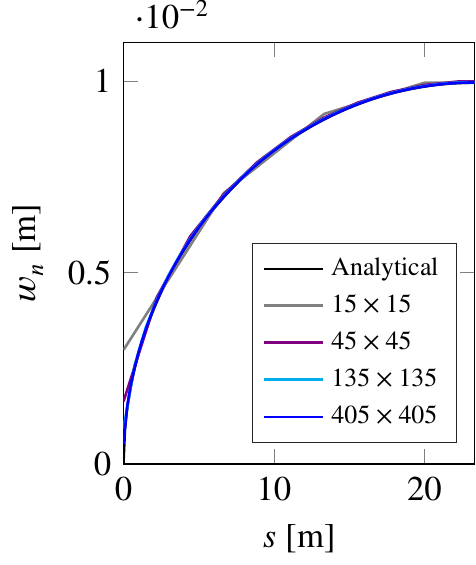}
		\caption{XFEM}
	\end{subfigure}	
	
	\caption{Test Case 1: Comparison of analytical and numerical aperture profiles for the different schemes applied to a horizontal fracture geometry.  Only half-profiles are shown due to symmetry.}
	\label{Fig:Case1H_solution}
\end{figure}

\begin{figure}[htbp]
	\newcommand{\panelsize}{0.3}
	
	\begin{subfigure}[t]{\panelsize\textwidth}  
		\centering
		\includegraphics[width=\linewidth]{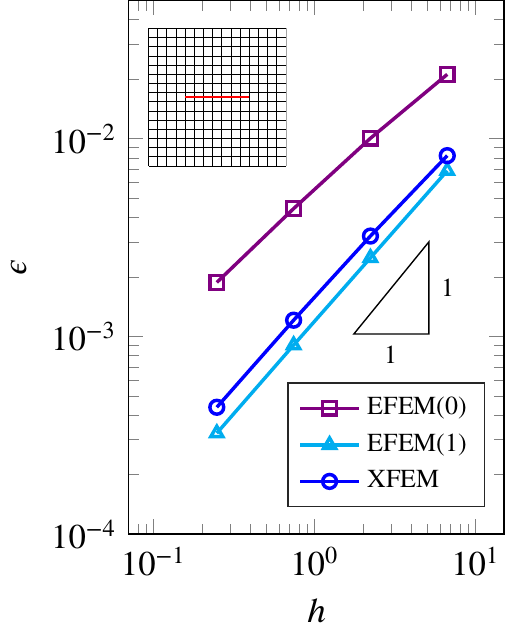}
		\caption{$\theta = 0\degree$}
	\end{subfigure}	
	\hfill	
	\begin{subfigure}[t]{\panelsize\textwidth}  
		\centering
		\includegraphics[width=\linewidth]{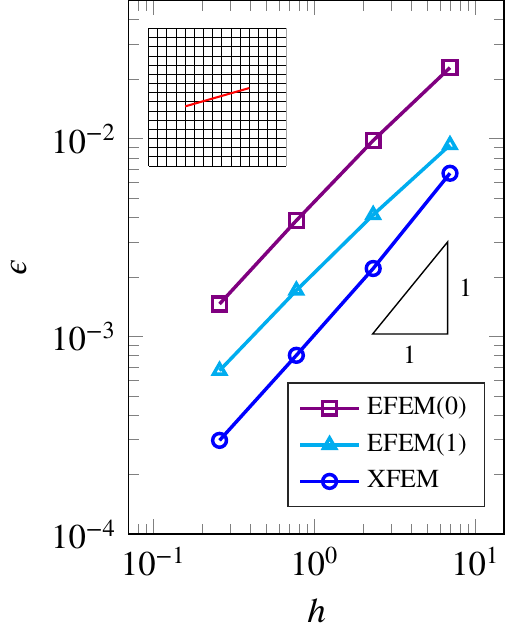}
		\caption{$\theta = 16\degree$}
	\end{subfigure}	
	\hfill	
	\begin{subfigure}[t]{\panelsize\textwidth}  
		\centering
		\includegraphics[width=\linewidth]{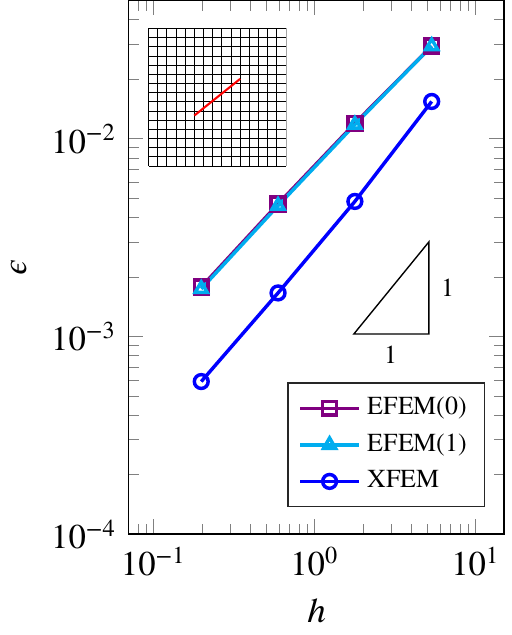}
		\caption{$\theta = 38\degree$}
	\end{subfigure}	  
	\caption{Test Case 1: Aperture error as a function of grid resolution for various orientation angles $\theta$ of the fracture. Base mesh and fracture location are shown at the top left corner in each panel. }
	\label{Fig:Case1_rot}
\end{figure}

Figure \ref{Fig:Case1H_solution} compares the analytical solution and the numerical one obtained for different grid resolutions, assuming a horizontal fracture. An aperture error is defined as 
\begin{equation}
\epsilon = \frac{\sqrt{\int\limits_0^{L} \left( w_n - w_n^{an} \right)^2 dL}}{\int\limits_0^{L} w_n^{an} dL} \,,
\end{equation} 
where $w_n^{an}$ is the analytical solution and $L$ is the total fracture length. Figure \ref{Fig:Case1_rot} presents the convergence behavior of the different schemes as a function of the rotation angle, which leads to different element intersection geometries.  

All three schemes exhibit linear or slightly superlinear convergence behavior, with the higher order schemes having a lower error constant.  For the horizontal fracture, EFEM(1) and XFEM have nearly the same accuracy.  As the fracture rotates, however, a larger number of single-cut-node elements are encountered.  As a result, the EFEM(1) scheme approaches the EFEM(0) behavior.  In general, however, all three approaches provide a good approximation of the opening profile.

\subsection{Test Case 2: Single fracture under compression} 

\begin{figure}[htbp]
	\centering
	\begin{tabular}{ccccc}
		\textbf{Uniform slip} &
		\hspace{.2\linewidth} & \multicolumn{3}{c}{\textbf{Non-Uniform slip}} \\
		\begin{tikzpicture}[scale=0.5]	
		\draw[color = white, pattern = north east lines] (0, 0) rectangle (4, -0.75);
		\draw [thin, black] (0,0) grid + (4, 8);
		\draw[very thick, violet] (0, 1.4) to  (4, 5.6);
		\node[align=center] at (2, 8.6) { $u_y = -0.1$};
		\end{tikzpicture} &
		\hspace{.2\linewidth} &
		\begin{tikzpicture}[scale=0.5]	
		\draw[color = white, pattern = north east lines] (0, 0) rectangle (4, -0.75);
		\draw [thin, black] (0,0) grid + (4, 8);
		\draw[very thick, violet] (0, 1.4) to  (2.6, 4);
		\node[align=center] at (2, 8.6) { $u_y = -0.1$};
		\end{tikzpicture} & 
		\begin{tikzpicture}[scale=0.5]	
		\draw[color = white, pattern = north east lines] (0, 0) rectangle (4, -0.75);
		\draw [thin, black, step = 0.5] (0,0) grid + (4, 8);
		\draw[very thick, violet] (0, 1.4) to  (2.6, 4);
		\node[align=center] at (2, 8.6) { $u_y = -0.1$};
		\end{tikzpicture} & 
		\begin{tikzpicture}[scale=0.5]	
		\draw[color = white, pattern = north east lines] (0, 0) rectangle (4, -0.75);
		\draw [thin, black, step=0.25] (0,0) grid + (4, 8);
		\draw[very thick, violet] (0, 1.4) to  (2.6, 4);
		\node[align=center] at (2, 8.6) { $u_y = -0.1$};
		\end{tikzpicture} \\ 
		&
		\hspace{.2\linewidth} & Grid 1 & Grid 2 & Grid 3 
	\end{tabular}
	\caption{Test Case 2: Geometry and boundary conditions for the uniform and non-uniform slip models.}
	\label{Fig:Case2_Geometry}
\end{figure}
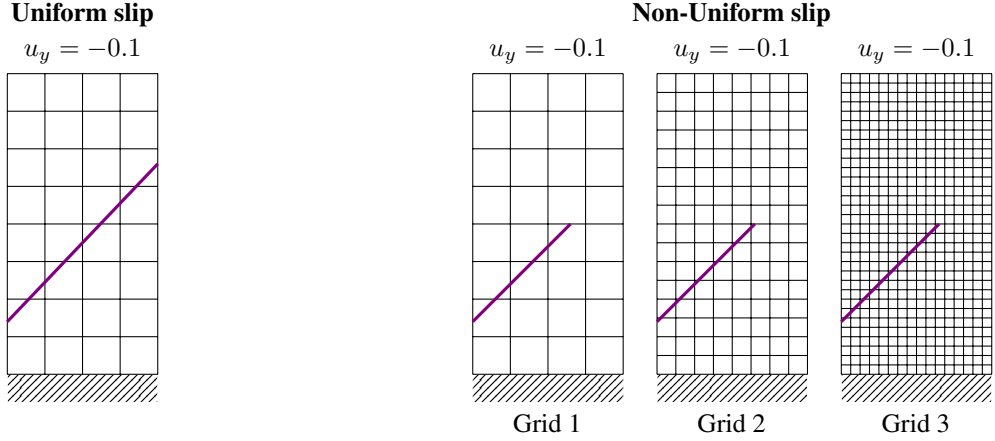

\begin{figure}[htbp]
	
	\newcommand{\panelsize}{0.45}
	
	\begin{subfigure}[t]{.82\linewidth}
		\begin{subfigure}[t]{\panelsize\linewidth}
			\includegraphics[width=.3\linewidth]{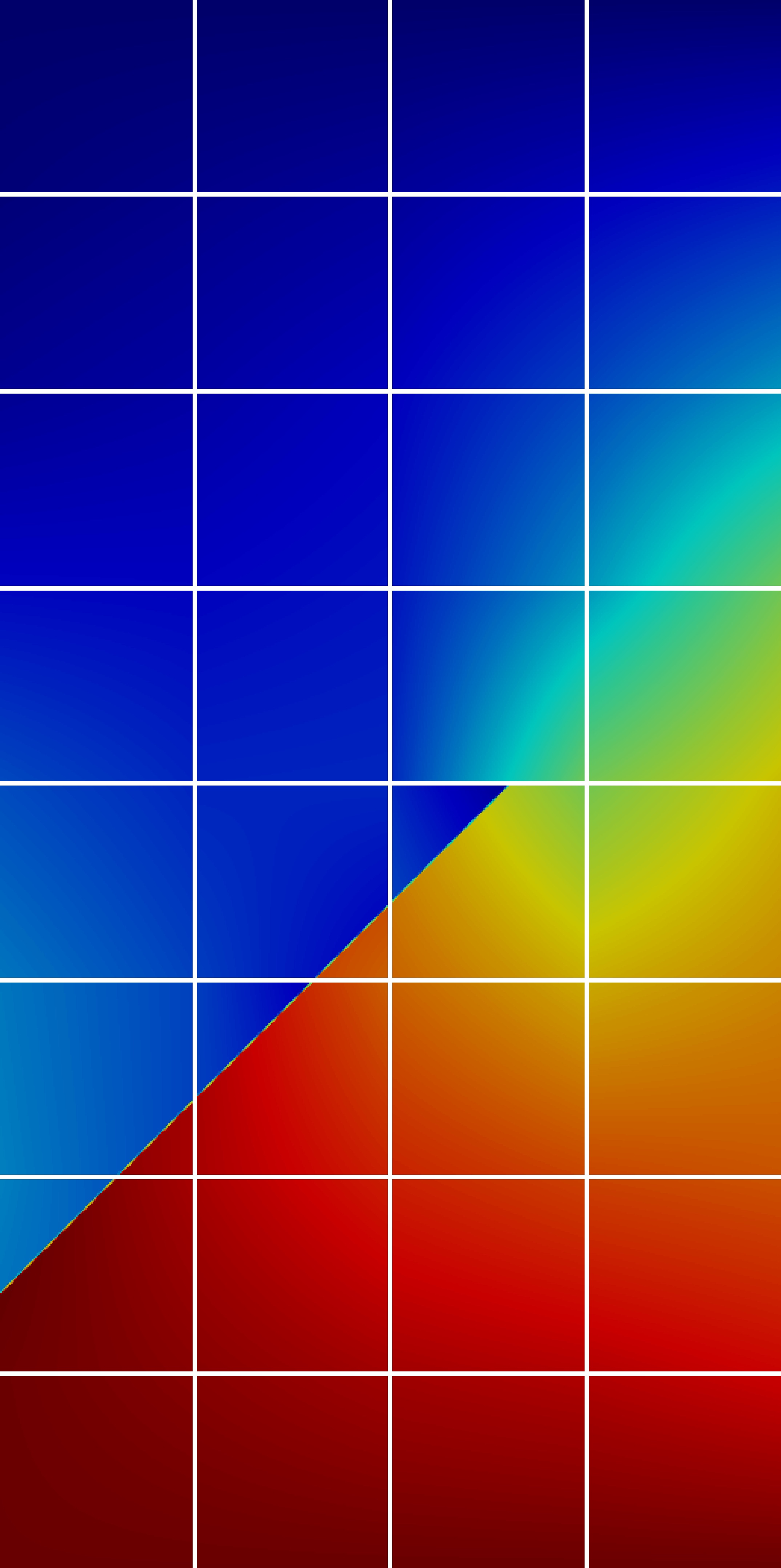}
			\hfill
			\includegraphics[width=.3\linewidth]{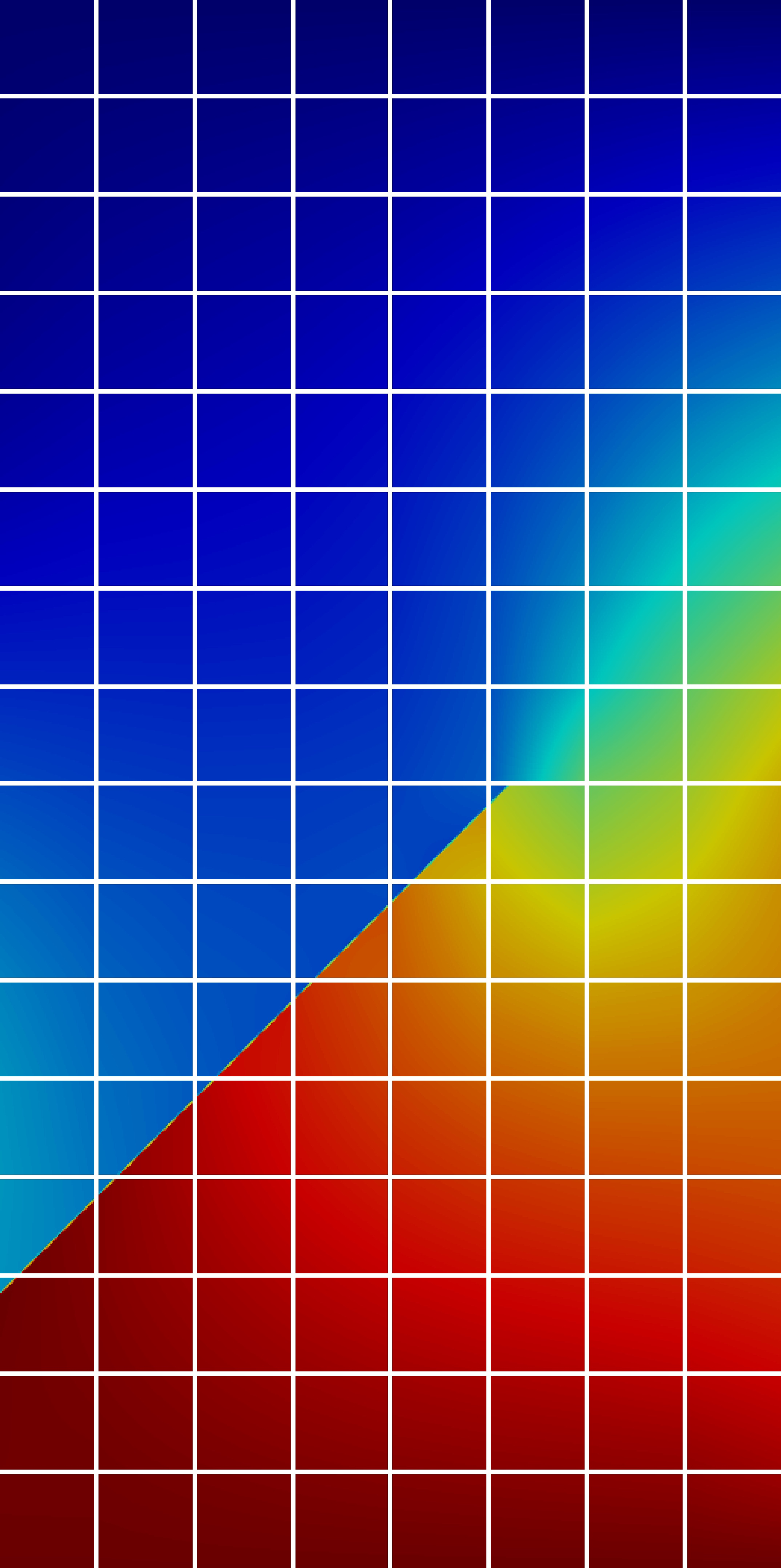}
			\hfill
			\includegraphics[width=.3\linewidth]{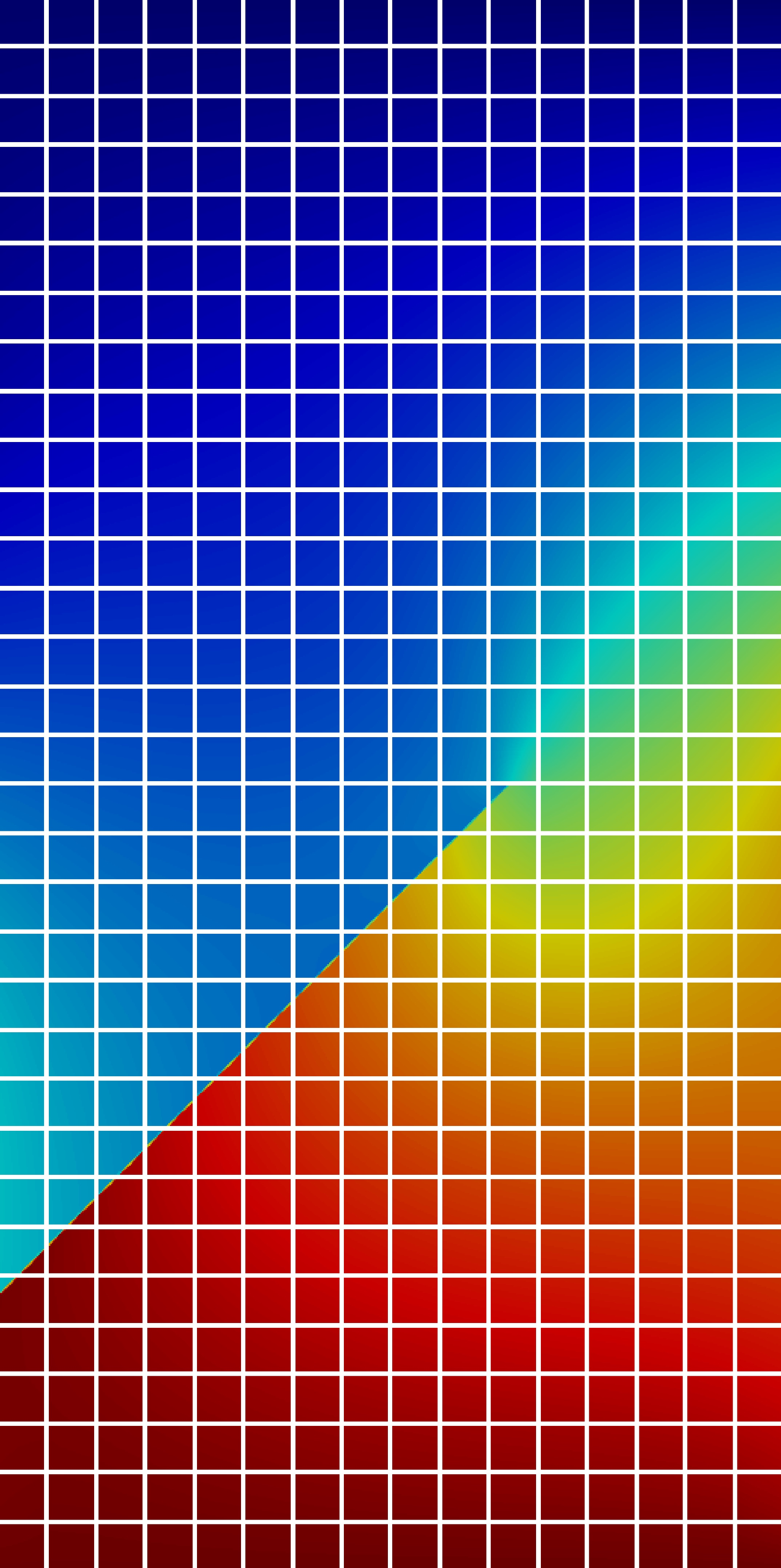}
			\caption{EFEM}
		\end{subfigure}	   
		\hfill
		\begin{subfigure}[t]{\panelsize\linewidth}
			\centering
			\includegraphics[width=.3\linewidth]{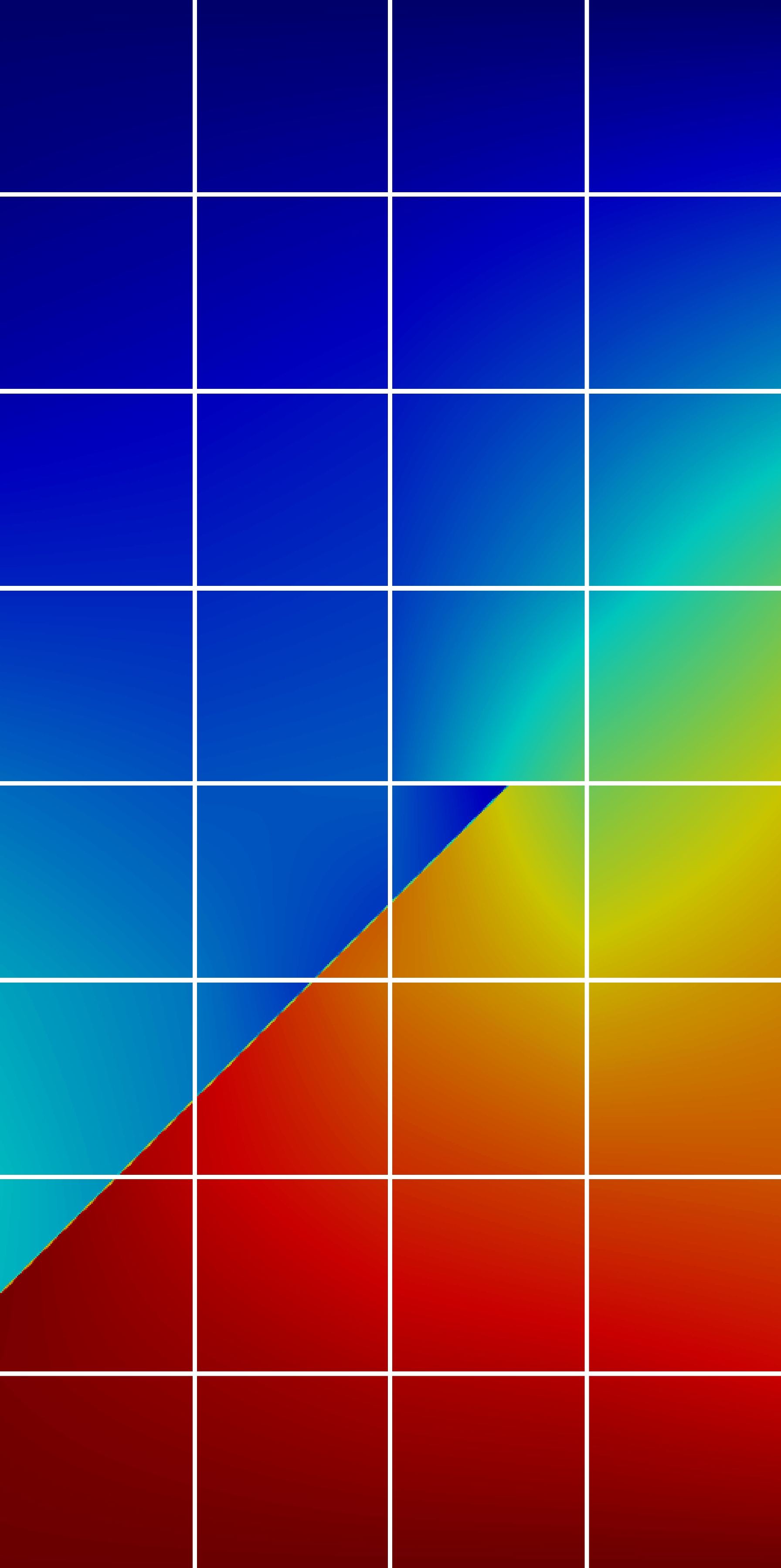}
			\hfill
			\includegraphics[width=.3\linewidth]{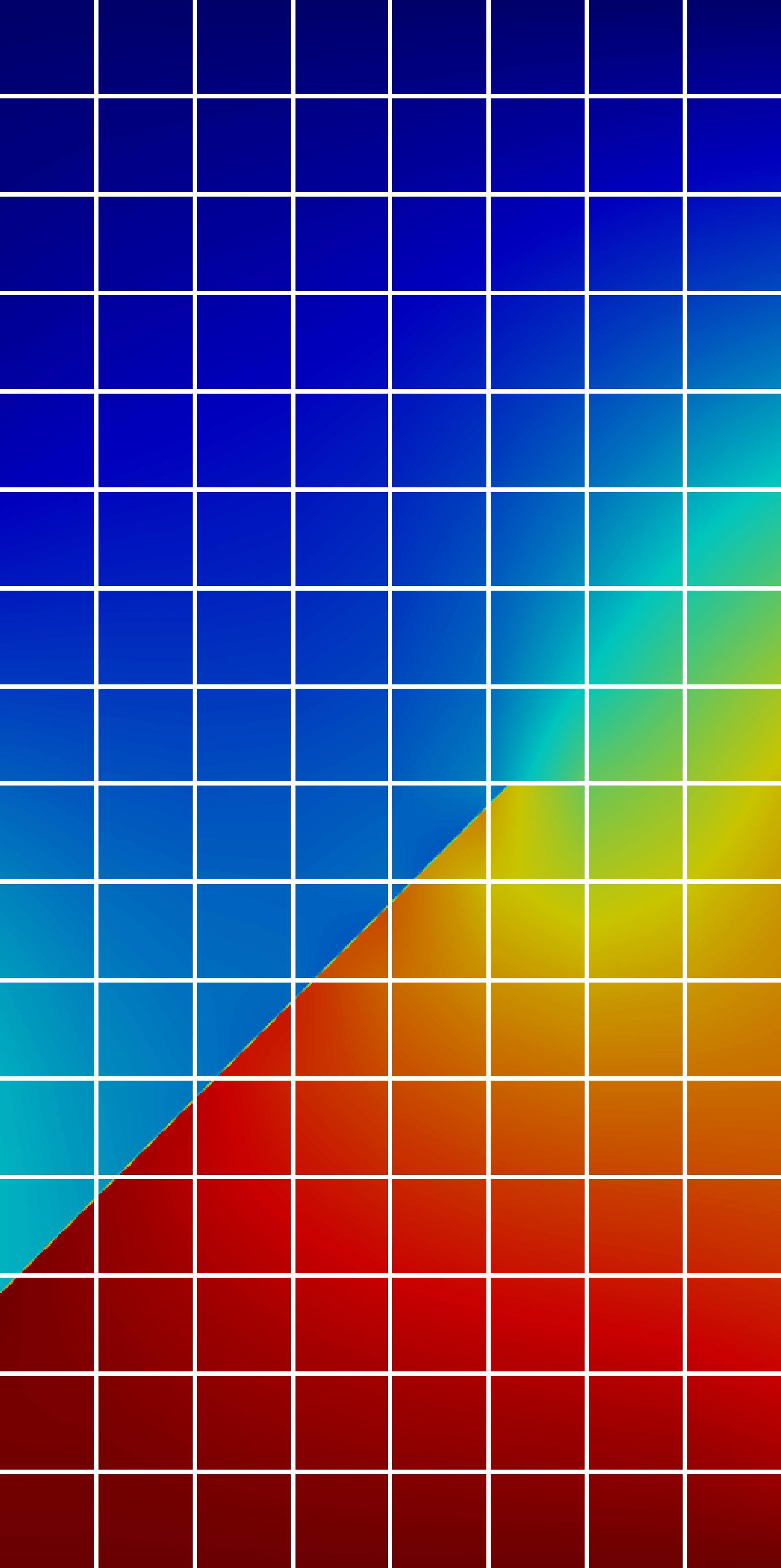}
			\hfill
			\includegraphics[width=.3\linewidth]{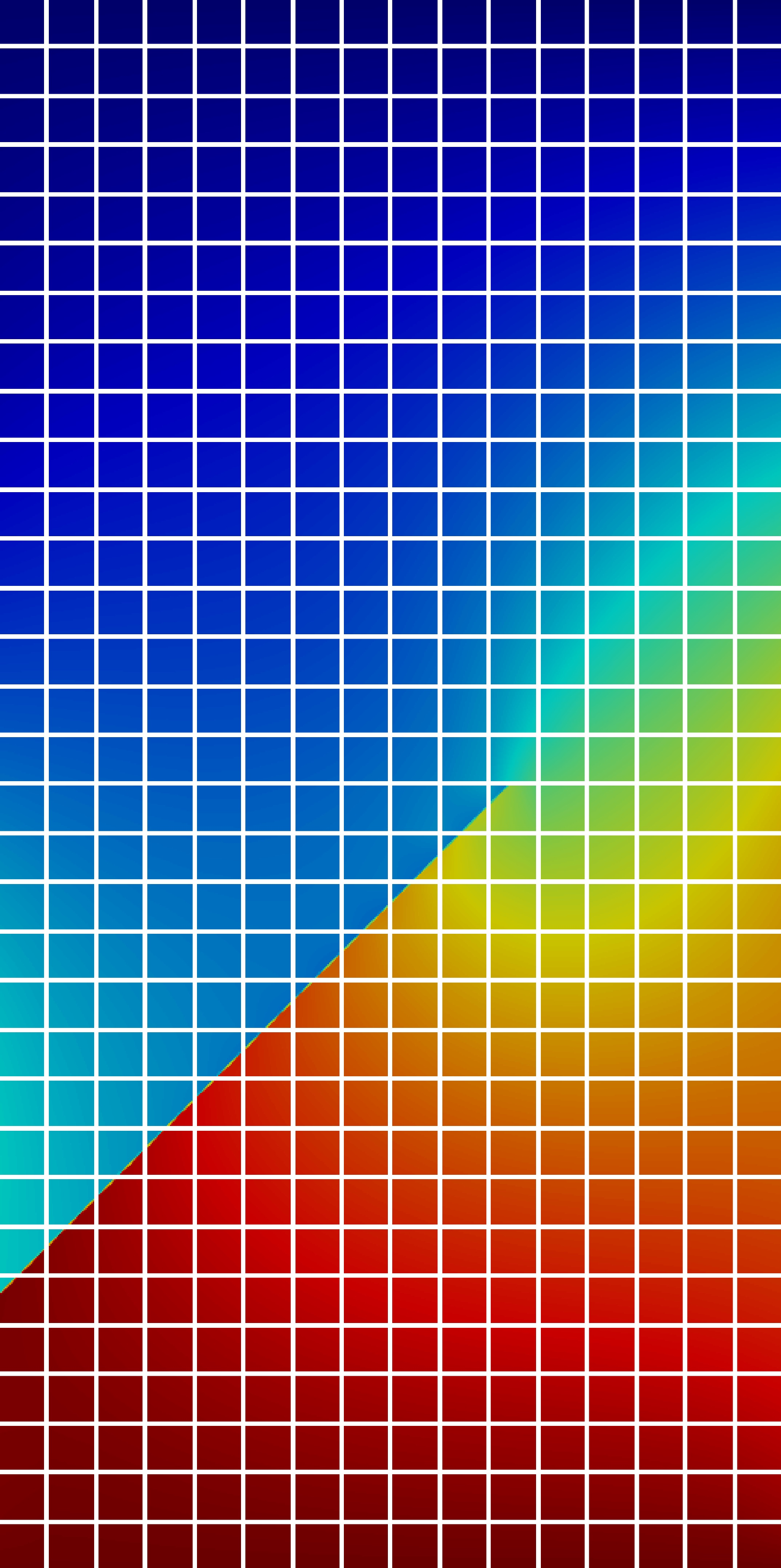}
			\caption{XFEM}
		\end{subfigure}	
	\end{subfigure}
	\hfill
	\begin{subfigure}[t]{.15\linewidth}
		\centering
		\begin{tikzpicture}
		\pgfplotscolorbardrawstandalone[
		colormap/jet,
		point meta min=-0.1,
		point meta max=0,
		parent axis height/.initial=2.75cm,
		colormap access=map,
		colorbar style={
			title=$u_y$ [m],
			scaled ticks=false,
			minor y tick num=1,
			/pgf/number format/precision=2,
			/pgf/number format/fixed,
			/pgf/number format/fixed zerofill
		}
		]
		\end{tikzpicture} 
	\end{subfigure}
	
	\caption{Test Case 2: Contours of vertical displacement for the non-uniform slip example.}
	\label{Fig:Case2_results}
\end{figure}

Here, we reproduce two examples originally presented in \cite{Borja2008}.  A $2 \, \text{m} \times 4 \, \text{m}$ elastic domain cut by a fracture inclined at $45 \degree$ is considered. In the first scenario, the fracture cuts through the entire domain, with end points at $\vec{x}_1=(0, 0.7)$ and $\vec{x}_2=(4, 2.7)$. Then, a shorter fracture is considered, which only cuts part of the domain, having end points at $\vec{x}_1=(0, 0.7)$ and $\vec{x}_2=(1.3, 2)$. The geometry for both cases is shown in Figure~\ref{Fig:Case2_Geometry}.  A displacement $u_y = - 0.1 \, \text{m}$ is imposed at the top surface while the bottom is fixed.  Given that all elements in this test are single-node-cut, we only consider the EFEM(0) and the XFEM schemes.  For the first geometry, a uniform slip, $w_t = 0.1 \sqrt{2}$ is obtained with a relative error of the order of $10^{-10}$. For the non-uniform slip case, contour plots of vertical displacement are presented in Figure \ref{Fig:Case2_results}. Both solutions are in good agreement with the original reference.

\subsection{Test Case 3: Injection in a saturated porous medium with an inclined fracture}

\begin{figure}[htbp]
	\centering
	\begin{subfigure}[b]{.4\linewidth}
		\centering
		\includegraphics[width=.9\linewidth]{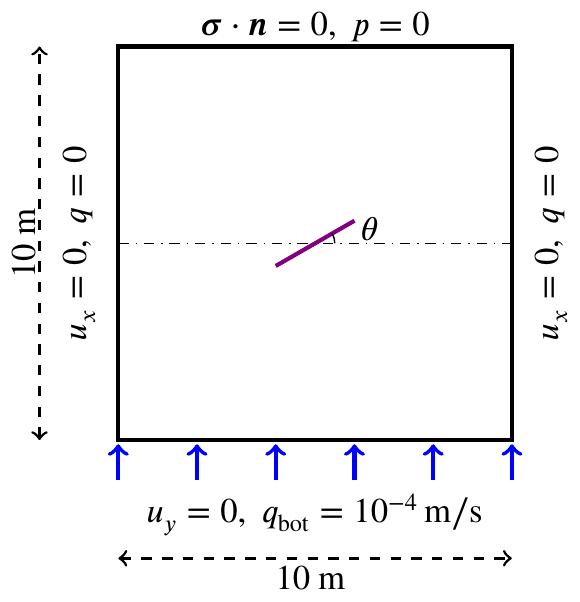}
		\caption{}
		\label{Fig:Case3_Geometry}
	\end{subfigure}
	\hfill
	\begin{subfigure}[b]{.5\linewidth}
		\centering
		\includegraphics[width=.9\linewidth]{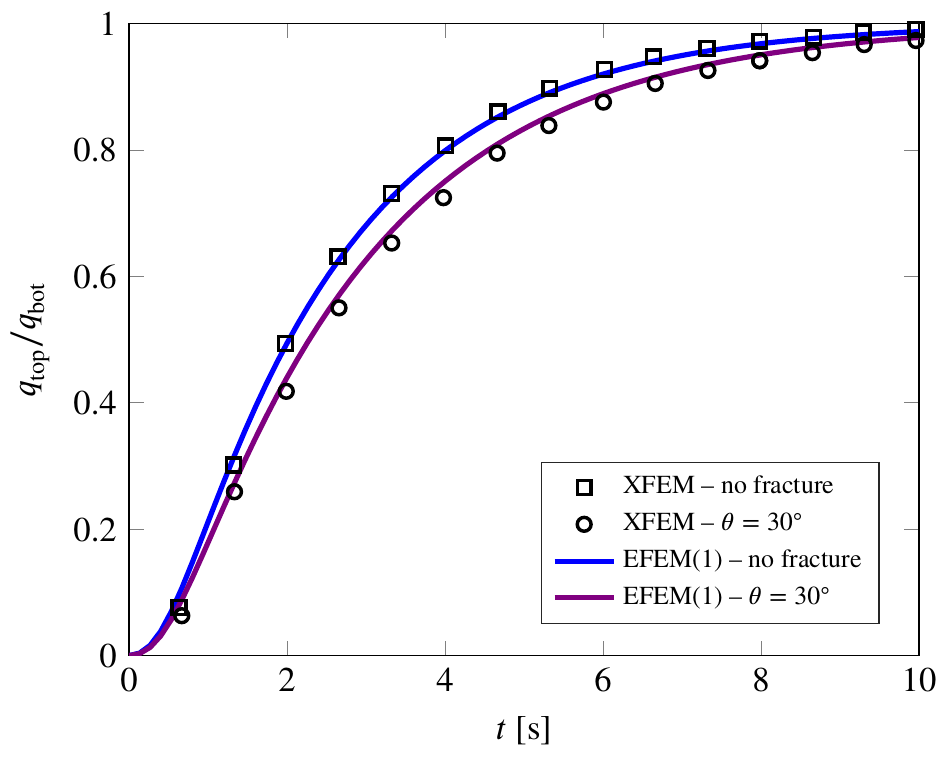}
		\caption{}
		\label{Fig:Case3_Results}
	\end{subfigure}	
	\caption{Test Case 3: (a) Geometry and boundary conditions. (b) Normalized flux at the top surface versus time, with and without a fracture inclined by $30\degree$ with respect to the horizontal direction. The black points were obtained using the XFEM-based method presented in \cite{Khoei2014} and provide an independent reference solution.}
\end{figure}

This test reproduces a benchmark example from \cite{Rethore2007,Khoei2014}.  A $10 \, \text{m} \times 10 \, \text{m}$ saturated reservoir with a $2$ m long inclined fracture positioned at the middle of the domain is considered (Figure \ref{Fig:Case3_Geometry}). The bottom surface is subject to a fluid influx at a constant rate, $q_\text{bot} = 10^{-4} \, {\text{m}}/{\text{s}}$, whereas a constant pressure boundary condition is applied at the top surface.  No flow is allowed at the two sides. Additionally, free displacement is allowed at the top surface, while all other surfaces are subject to zero displacement in the normal direction.  

The initial reservoir pressure is $p_\text{0} = 0 \, \text{MPa}$ and the injection process is run for 10 seconds. Three different configurations are considered: one in which no fracture is present and two configurations in which the fracture is inclined by $30\degree$ and $60\degree$ with respect to the horizontal axis. A $30 \times 30$ cartesian grid is imposed on the domain and the simulation is run using 75 time-steps for a total simulation time of 10 s. 

\begin{figure}[htbp]
	\newcommand{\panelsize}{0.25}
	
	\hfill
	\begin{subfigure}[t]{\panelsize\linewidth}
		\centering
		\includegraphics[width=\linewidth]{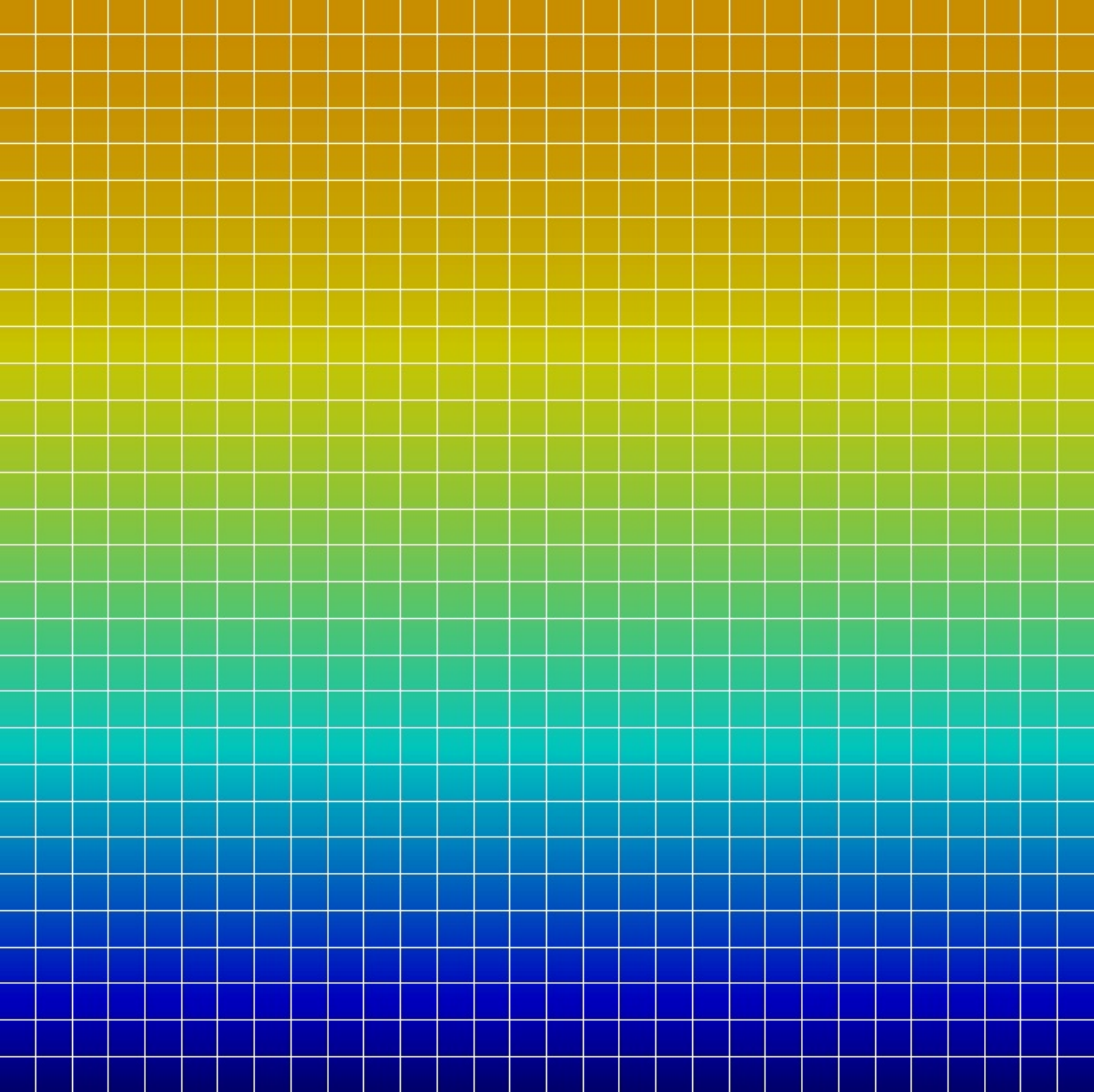}
		\caption{No fracture}
	\end{subfigure} 
	\hfill
	\begin{subfigure}[t]{\panelsize\linewidth}
		\centering
		\includegraphics[width=\linewidth]{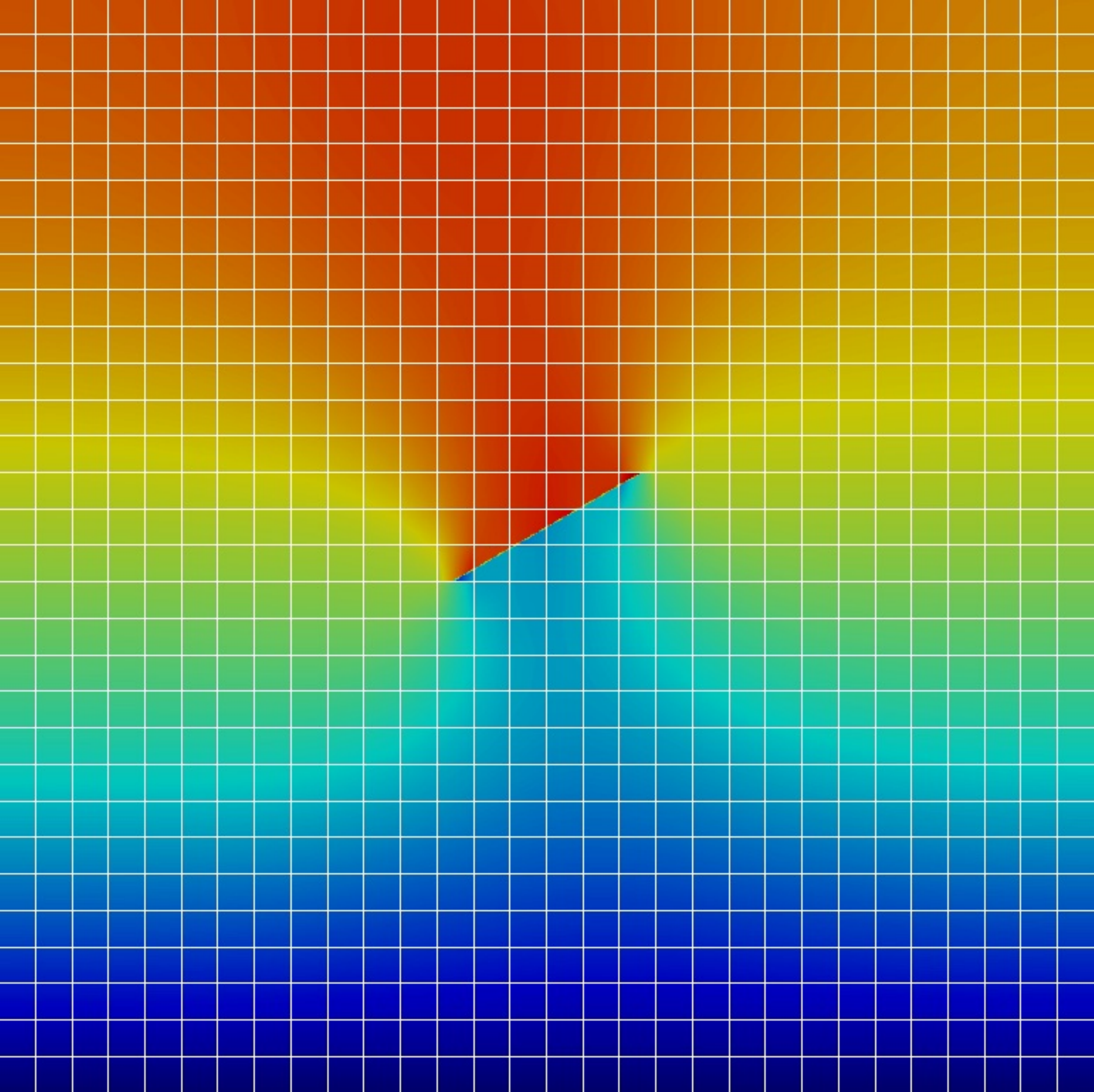}
		\caption{$\theta = 30\degree$}
	\end{subfigure}
	\hfill   
	\begin{subfigure}[t]{\panelsize\linewidth}
		\centering
		\includegraphics[width=\linewidth]{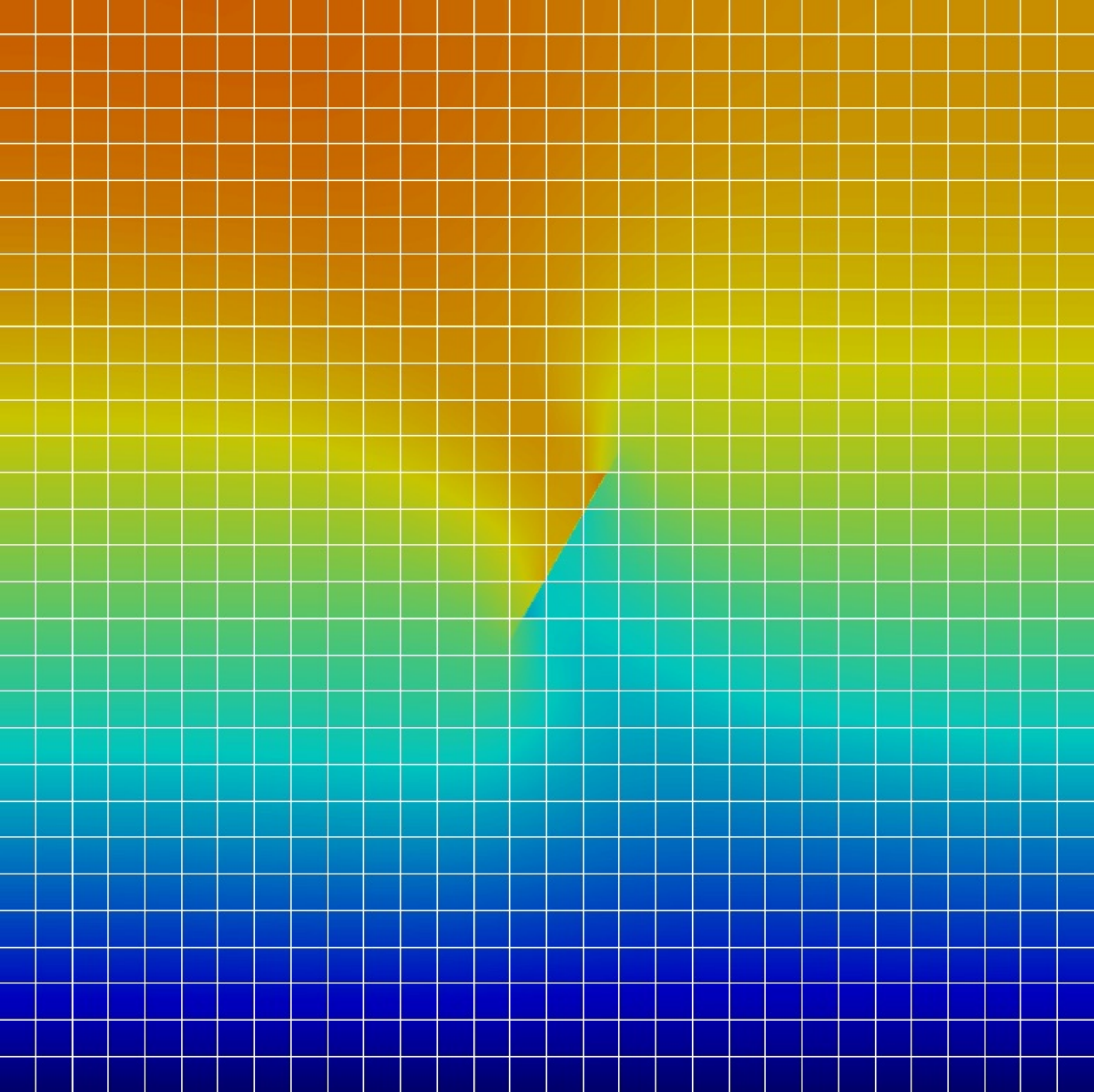}
		\caption{$\theta = 60 \degree$}
	\end{subfigure}
	\hfill  
	\begin{subfigure}[t]{.15\linewidth}
		\centering
		\begin{tikzpicture}
		\pgfplotscolorbardrawstandalone[
		colormap/jet,
		point meta min=0,
		point meta max=3.723e-4,
		parent axis height/.initial=3cm,
		colormap access=map,
		colorbar style={
			title=$u_y$ [m],
			scaled ticks=false,
			minor y tick num=1,
			/pgf/number format/precision=2
		}
		]
		\end{tikzpicture} 
	\end{subfigure}	
	\hfill\null
	\caption{Test Case 3: Contour plots of vertical displacement for three scenarios.}
	\label{Fig:Case3_Disp}
\end{figure}

The example is run with the EFEM(1) method, and compared to the XFEM-based results presented by the original authors of the benchmark \cite{Rethore2007,Khoei2014}.  We note that their method employs a different discretization strategy for the flow equations.  Figure \ref{Fig:Case3_Results} plots the outgoing flux at the top surface $q_\text{top}$ as a function of time, normalized by the injection rate $q_\text{bot}$. The fracture, due to its storage capacity, delays the time that it takes for the flux at the top surface to reach steady state. Figure \ref{Fig:Case3_Disp} plots the vertical displacement at the end of the simulation. As expected, the vertical displacement is larger for the smallest value of the angle $\theta$. 

\subsection{Test Case 4: Primary depletion and water injection in a fractured reservoir}
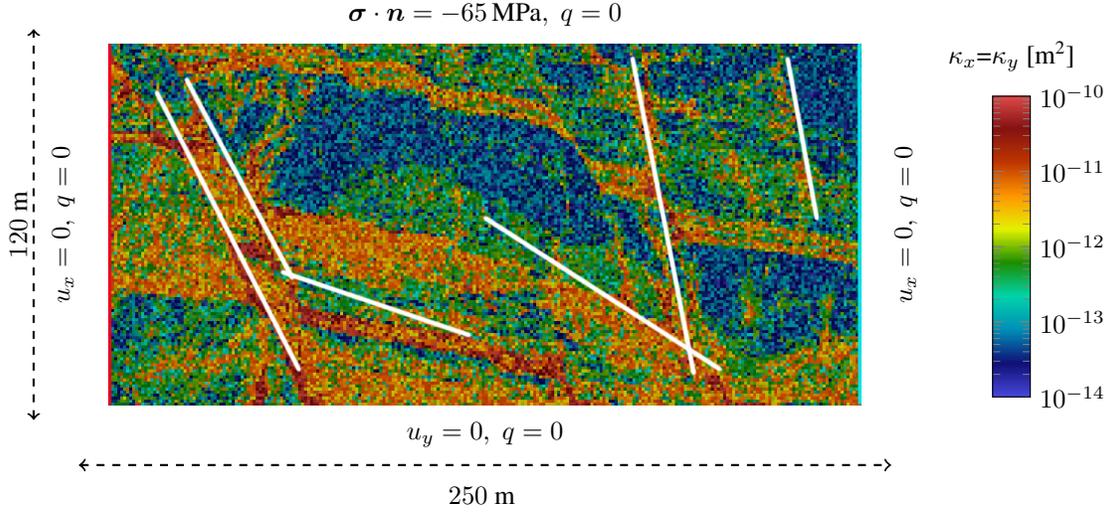
\begin{figure}[htbp]
	\centering
	\begin{tikzpicture}[scale=0.4]
	
	\node at (0, 0) {\includegraphics[width =  0.606\textwidth]{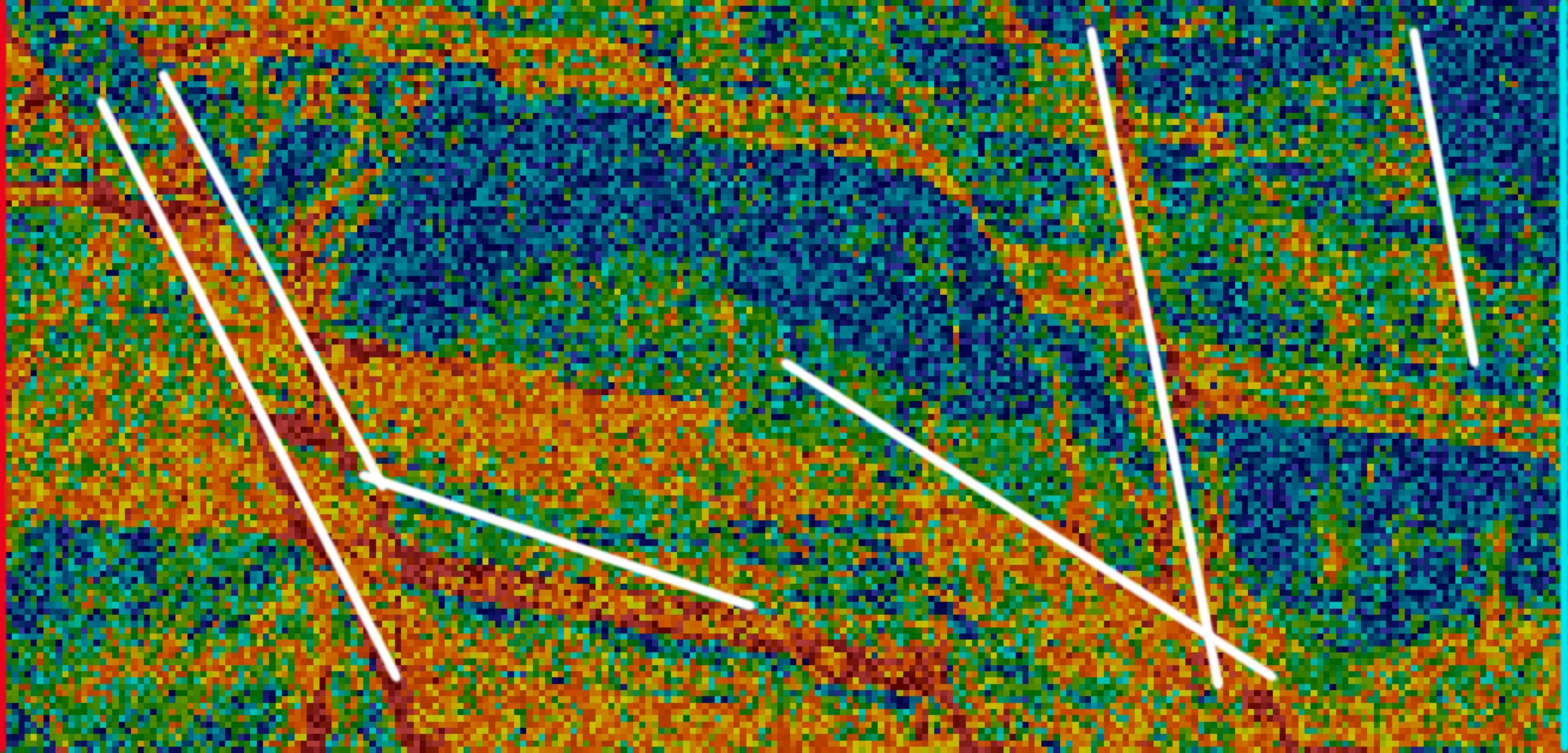}};
	\node at (18, 0) {	    \pgfplotscolorbardrawstandalone[
		colormap={example}{
			samples of colormap=(10 of hotDesaturated)
		},
		point meta min=1.003e-14,
		point meta max=9.984e-11,
		parent axis height/.initial=4cm,
		colormap access=map,
		colorbar style={
			title={$\kappa_x$=$\kappa_y$ [m\textsuperscript{2}]},
			scaled ticks=false,
			ymode=log
		}
		]};
	
	\draw[thick, dashed, <->] (-15, -6.5) -- (-15, 6.5);
	\draw[thick, dashed, <->] (-13.5, -8) -- (13.5, -8);
	\node[align=center, rotate = 90] at (-15.5, 0) {$120$ m};
	\node[align=center] at (0, -9) {$250$ m };
	
	\node[align=center] at (0,-7) {$u_y = 0, \; q = 0$};
	\node[align=center] at (0, 7) {$\tensorTwo{\sigma} \cdot \vec{n} =-65 \, \text{MPa}, \; q = 0$};
	\node[align=center, rotate = 90] at (-14,0) {$u_x = 0,\; q = 0$};
	\node[align=center, rotate = 90] at ( 14,0) {$u_x = 0, \; q = 0$};

	\end{tikzpicture}
	\caption{Test case 4: Geometry, boundary conditions, and permeability map in logarithmic scale. Fractures (white), producer 1 (red) and producer 2 (cyan) are also shown.}
	\label{Fig:Case4_GeometryPerm}
\end{figure}

A vertical section of a heterogeneous fractured reservoir is considered. The reservoir is $120\, \text{m}$ high and has a length of $250\,\text{m}$. The permeability field and the boundary conditions are shown in Figure \ref{Fig:Case4_GeometryPerm}. The domain is discretized with a $250\times120$ cartesian grid. Six large fractures are embedded in the domain, and two wells are introduced at the left and right-hand sides of the reservoir. Here, the Peaceman well model is employed for the wells \cite{Peaceman1978}.  Both wells are pressure-constrained with a pressure, $p_{\text{prod}} = 5 \, \text{MPa}$. The coordinates of the end points of each fracture and the cells perforated by each well are provided in tables \ref{Tab:FracturesLocation} and \ref{Tab:WellsLocation}. The reservoir initially contains two compressible phases with saturations $S_w = 0.1$ and $S_o = 0.9$, respectively. The initial reservoir pressure is uniform and equal to $p_{\text{init}} = 20 \, \text{MPa}$, and all fractures are closed. The initial vertical and horizontal effective stresses are $\sigma_V' = 45 \, \text{MPa}$ and $\sigma_H' = 0.33 \sigma_V'$. 

Two scenarios are considered. In the first, the reservoir is subject to primary depletion via fluid production through the two wells. In the second scenario, water is injected into fracture no. 5, increasing the average reservoir pressure and forcing some of the fractures to open.  Given the comparable performance of the various discretizations, only results for the EFEM(1) scheme are presented.

\subsubsection{Scenario 1: Primary depletion}
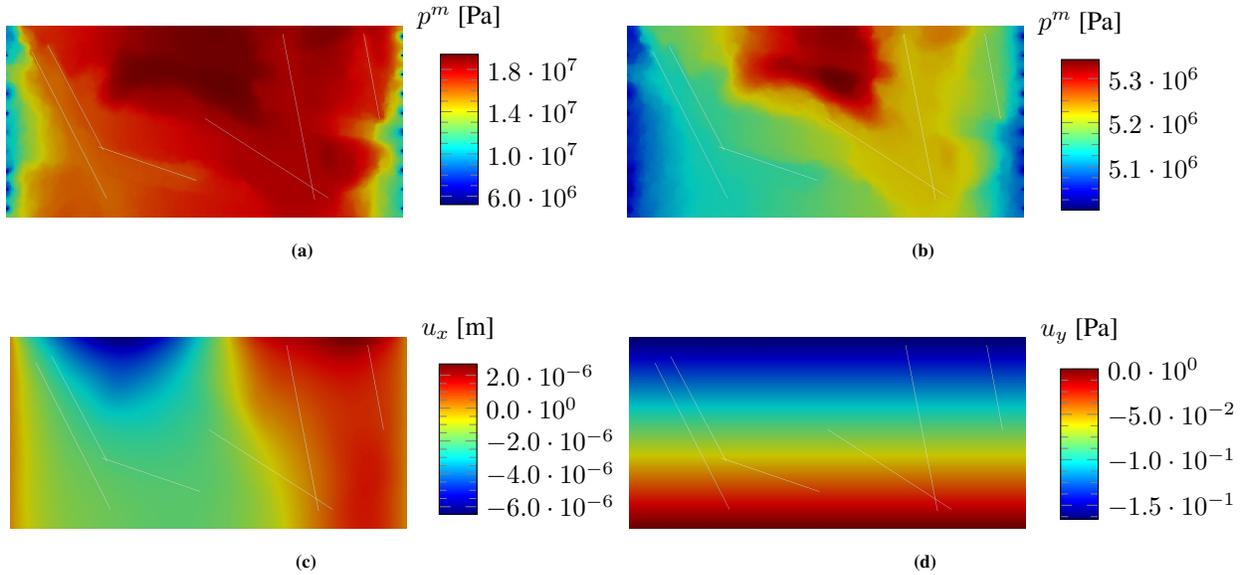
\begin{figure}[htbp]
	
	\centering
	
	\newcommand{\panelsize}{0.49}

	\begin{subfigure}[b]{\panelsize\linewidth}	
		\flushleft    
		
		\begin{tikzpicture}
		\node at (0, 0) {\includegraphics[width=.65\textwidth]{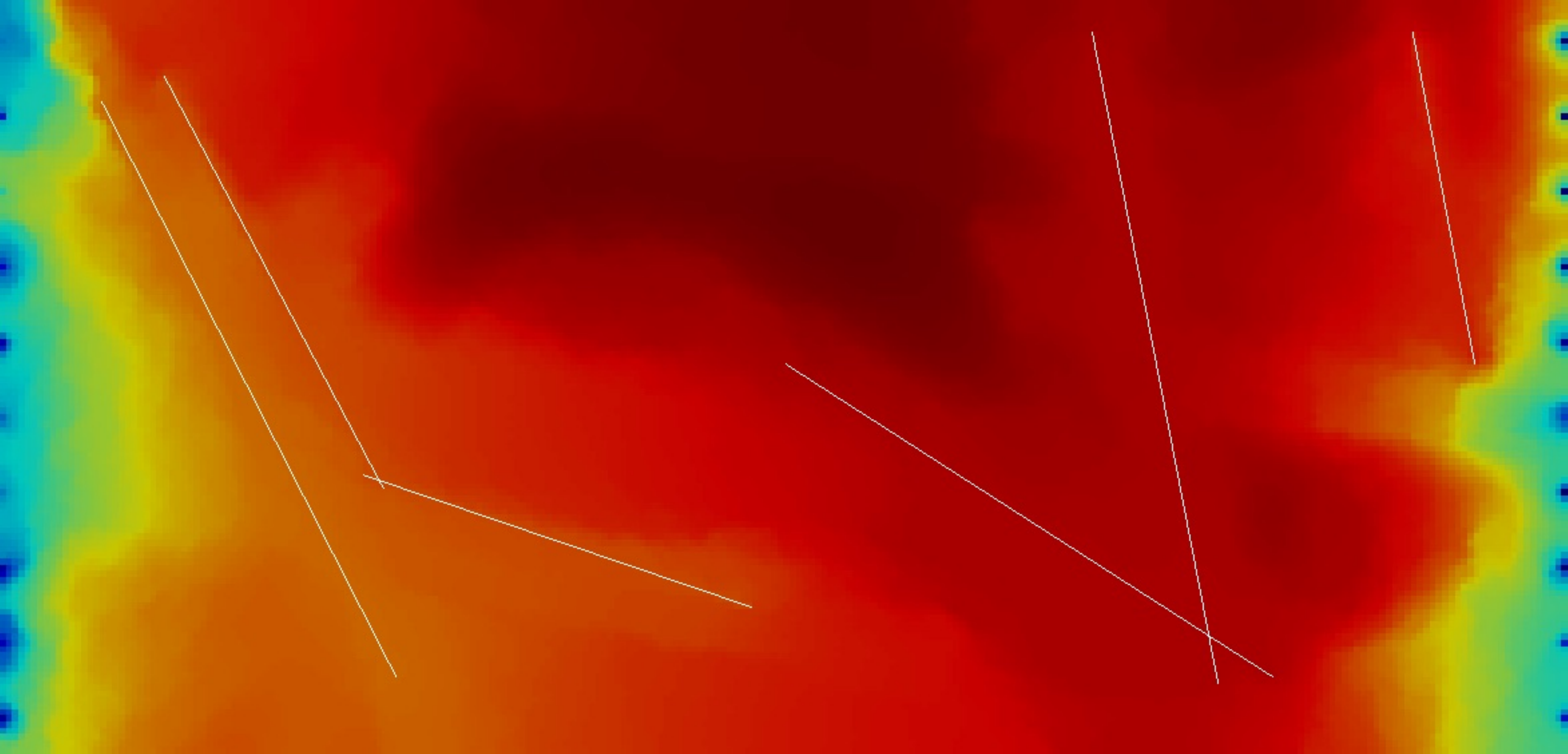}};
		\node[anchor=west] at (.32\linewidth, 2mm) {	      
			
			\pgfplotscolorbardrawstandalone[
			colormap/jet,
			point meta min=5.204e6,
			point meta max=1.944e7,
			parent axis height/.initial=2.cm,
			colormap access=map,
			colorbar style={
				title=$p^m$ [Pa],
				scaled ticks=false,
				ytick={0.6e7,1.0e7,1.4e7,1.8e7},
				minor y tick num = 3,
				/pgf/number format/sci,
				/pgf/number format/sci zerofill,
				/pgf/number format/precision=1
			}
			]};
		\end{tikzpicture}
		\caption{}
	\end{subfigure}	
	\hfill
	\begin{subfigure}[b]{\panelsize\linewidth}	
		\flushleft    
		
		\begin{tikzpicture}
		\node at (0, 0) {\includegraphics[width=.65\textwidth]{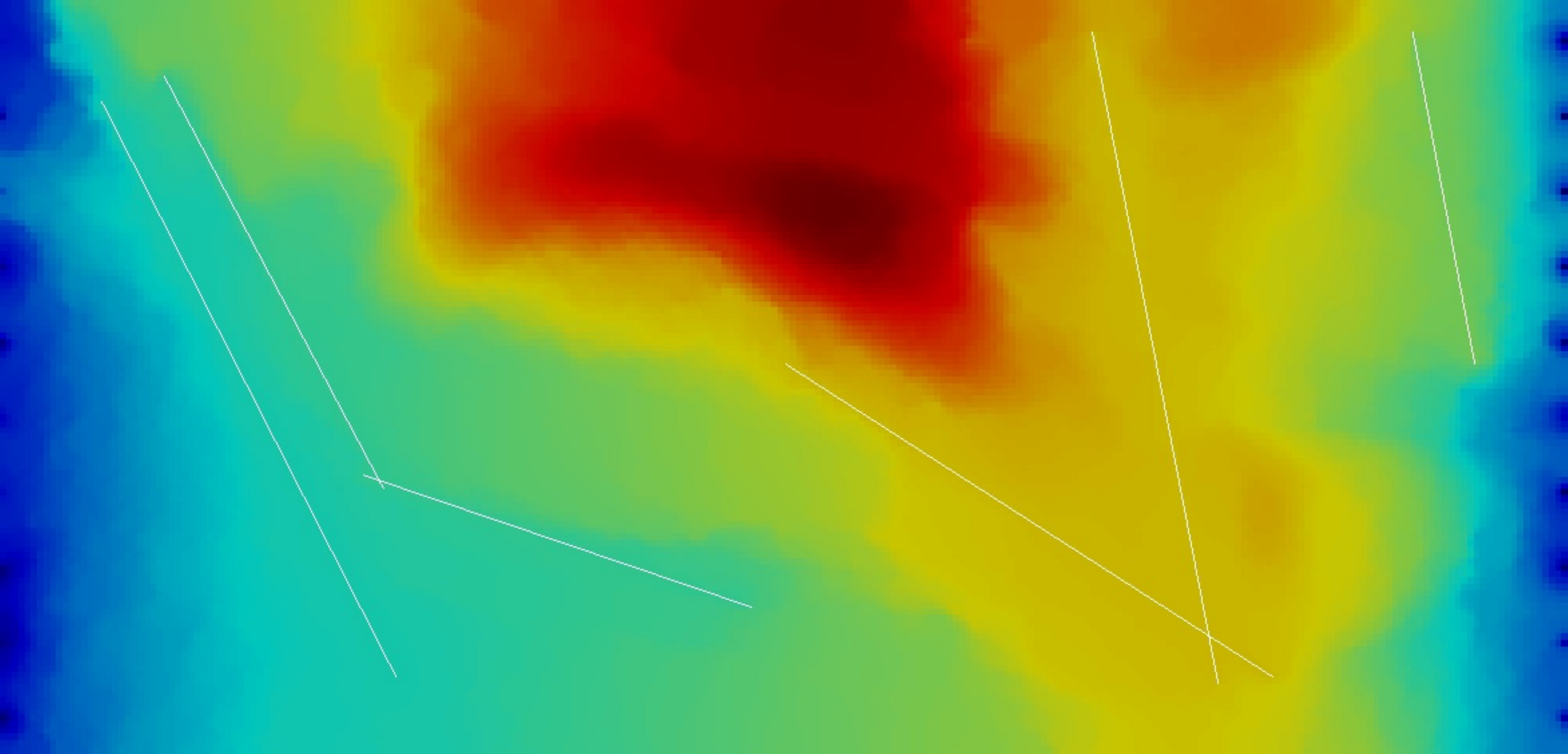}};
		\node[anchor=west] at (.32\linewidth, 2mm) {	  
			
			\pgfplotscolorbardrawstandalone[
			colormap/jet,
			point meta min=5.002e6,
			point meta max=5.343e6,
			parent axis height/.initial=2.cm,
			colormap access=map,
			colorbar style={
				title=$p^m$ [Pa],
				scaled ticks=false,
				minor y tick num = 3,
				/pgf/number format/sci,
				/pgf/number format/sci zerofill,
				/pgf/number format/precision=1
			}
			]};
		\end{tikzpicture}
		\caption{}
	\end{subfigure}	  
	\hfill\null
	
	\bigskip
	
	\hfill
	\begin{subfigure}[b]{\panelsize\linewidth}	
		\flushleft    
		
		\begin{tikzpicture}
		\node at (0, 0) {\includegraphics[width=.65\textwidth]{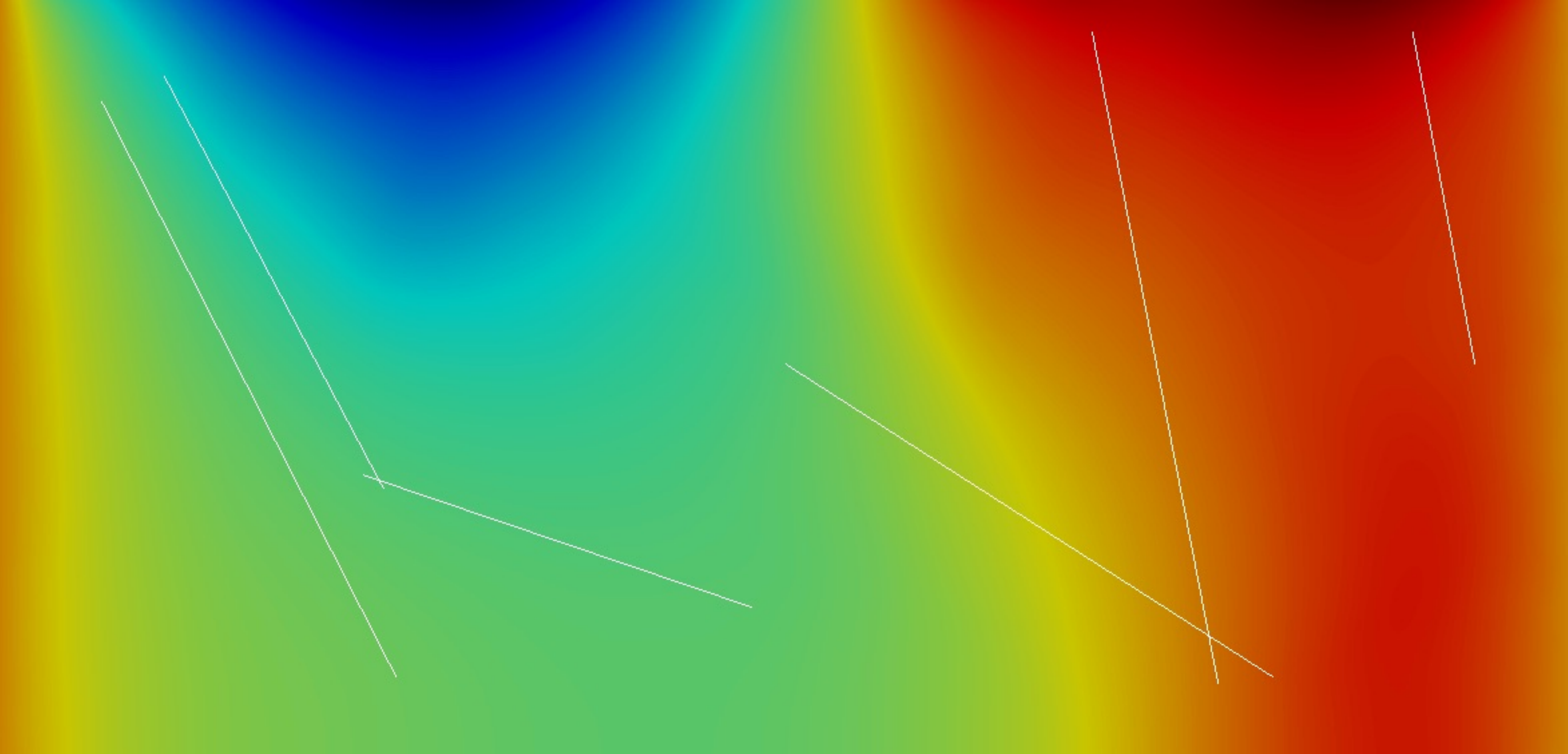}};
		\node[anchor=west] at (.32\linewidth, 2mm) {	  
			
			\pgfplotscolorbardrawstandalone[
			colormap/jet,
			point meta min=-6.461e-6,
			point meta max=2.688e-6,
			parent axis height/.initial=2.cm,
			colormap access=map,
			colorbar style={
				title=$u_x$ [m],
				scaled ticks=false,
				minor y tick num = 3,
				/pgf/number format/sci,
				/pgf/number format/sci zerofill,
				/pgf/number format/precision=1
			}
			]};
		\end{tikzpicture}
		\caption{}
	\end{subfigure}	
	\hfill
	\begin{subfigure}[b]{\panelsize\linewidth}	
		\flushleft     
		
		\begin{tikzpicture}
		\node at (0, 0) {\includegraphics[width=.65\textwidth]{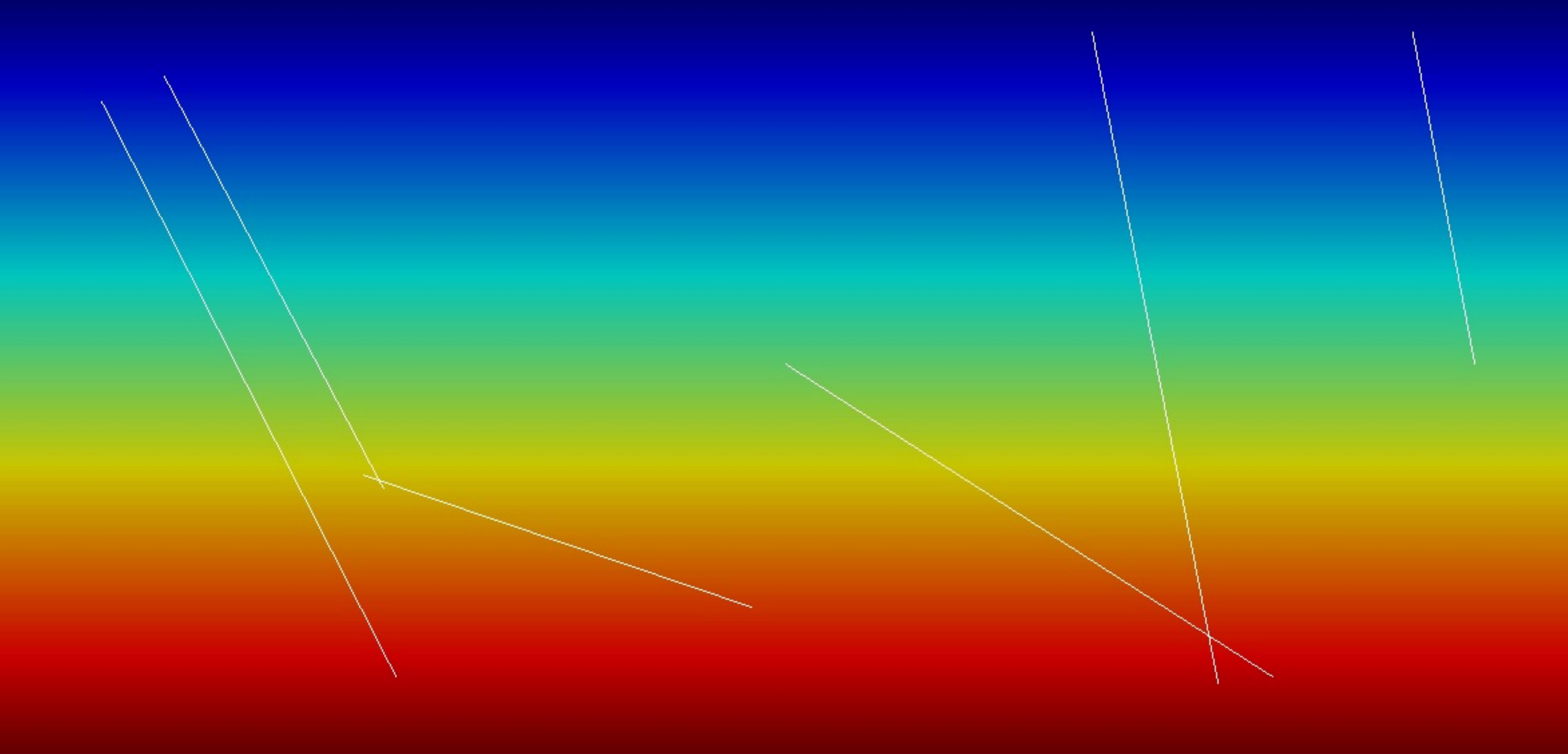}};
		\node[anchor=west] at (.32\linewidth, 2mm) {	  
			
			\pgfplotscolorbardrawstandalone[
			colormap/jet,
			point meta min=-1.653e-1,
			point meta max=0,
			parent axis height/.initial=2.cm,
			colormap access=map,
			colorbar style={
				title=$u_y$ [Pa],
				scaled ticks=false,
				minor y tick num = 3,
				/pgf/number format/sci,
				/pgf/number format/sci zerofill,
				/pgf/number format/precision=1
			}
			]};
		\end{tikzpicture}
		\caption{}
	\end{subfigure}	  
	\hfill\null 
	
	\caption{Test Case 4, Scenario 1: Contour plot of the reservoir pressure at time-steps 1 (a) and 8 (b) and of the displacement in $x$ (c) and in $y$ (d) at the end of the simulation.}
	\label{Fig:Case4Pressure}
\end{figure}

The reservoir pressure is decreased by extracting the two fluid phases from the production wells. The total simulation time is equal to $T_\text{final} = 1.33 \, \text{days}$ using 10 time-steps. The initial time-step size is 0.01 days and it is gradually increased throughout the simulation based on the Newton convergence rate. Figure \ref{Fig:Case4Pressure} shows the pressure map at time-steps 1 (0.01 days) and 8 (0.37 days). Figure  presents contour plots of the x and y components of the displacement field at the end of the simulation. As expected, the decrease in reservoir pressure results in subsidence at the top surface. 

\subsubsection{Scenario 2: Water injection in a fracture}
\begin{figure}[htbp]
	\centering
	
	\newcommand{\panelsize}{0.49}

	\begin{subfigure}[b]{\panelsize\linewidth}	
		\flushleft    
		
		\begin{tikzpicture}
		\node at (0, 0) {\includegraphics[width=.65\textwidth]{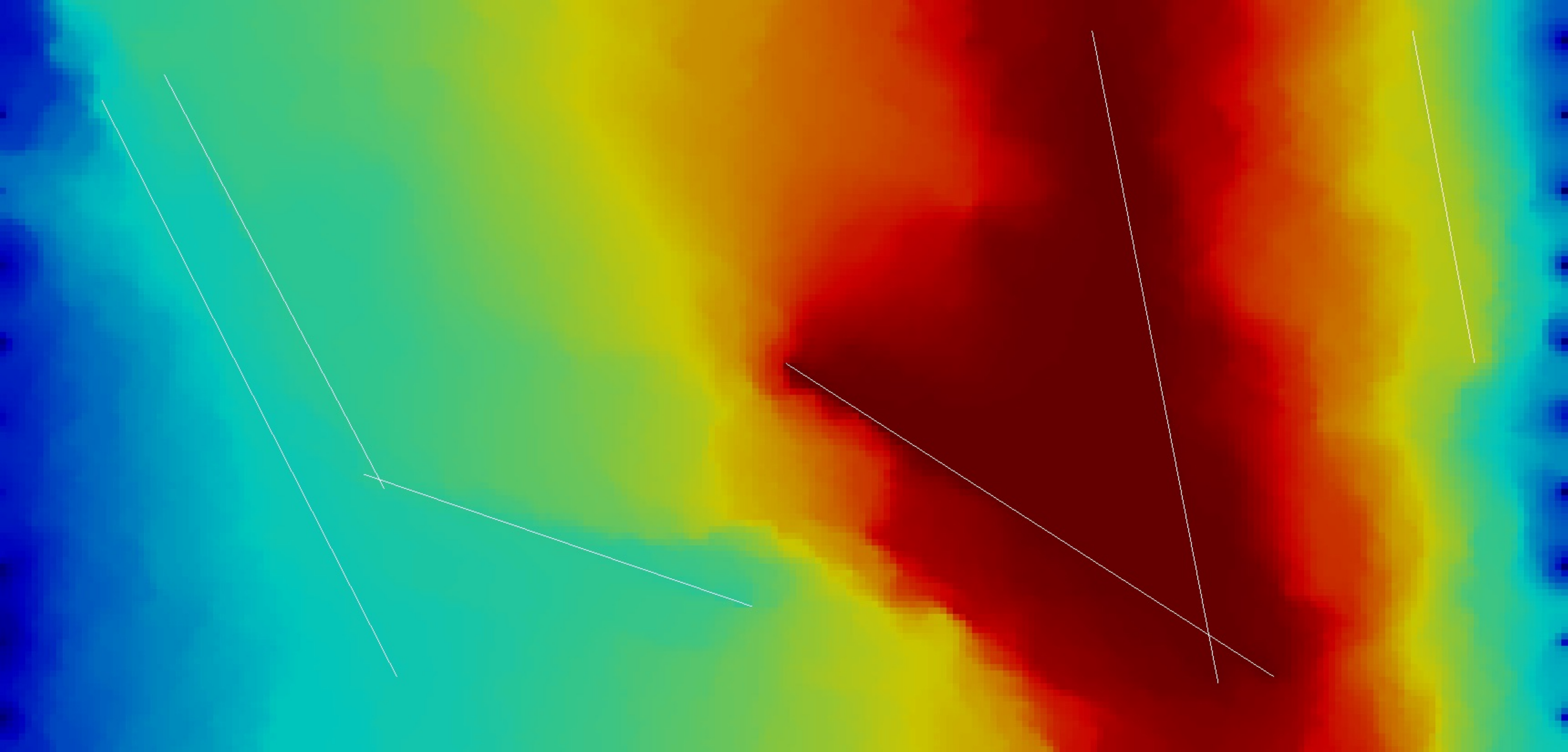}};
		\node[anchor=west] at (.32\linewidth, 2mm) {	  
			
			\pgfplotscolorbardrawstandalone[
			colormap/jet,
			point meta min=5.365e6,
			point meta max=7.500e7,
			parent axis height/.initial=2.cm,
			colormap access=map,
			colorbar style={
				title=$p^m$ [Pa],
				scaled ticks=false,
				minor y tick num = 3,
				/pgf/number format/sci,
				/pgf/number format/sci zerofill,
				/pgf/number format/precision=1
			}
			]};
		\end{tikzpicture}
		\caption{}
	\end{subfigure}	
	\hfill
	\begin{subfigure}[b]{\panelsize\linewidth}	
		\flushleft    
		
		\begin{tikzpicture}
		\node at (0, 0) {\includegraphics[width=.65\textwidth]{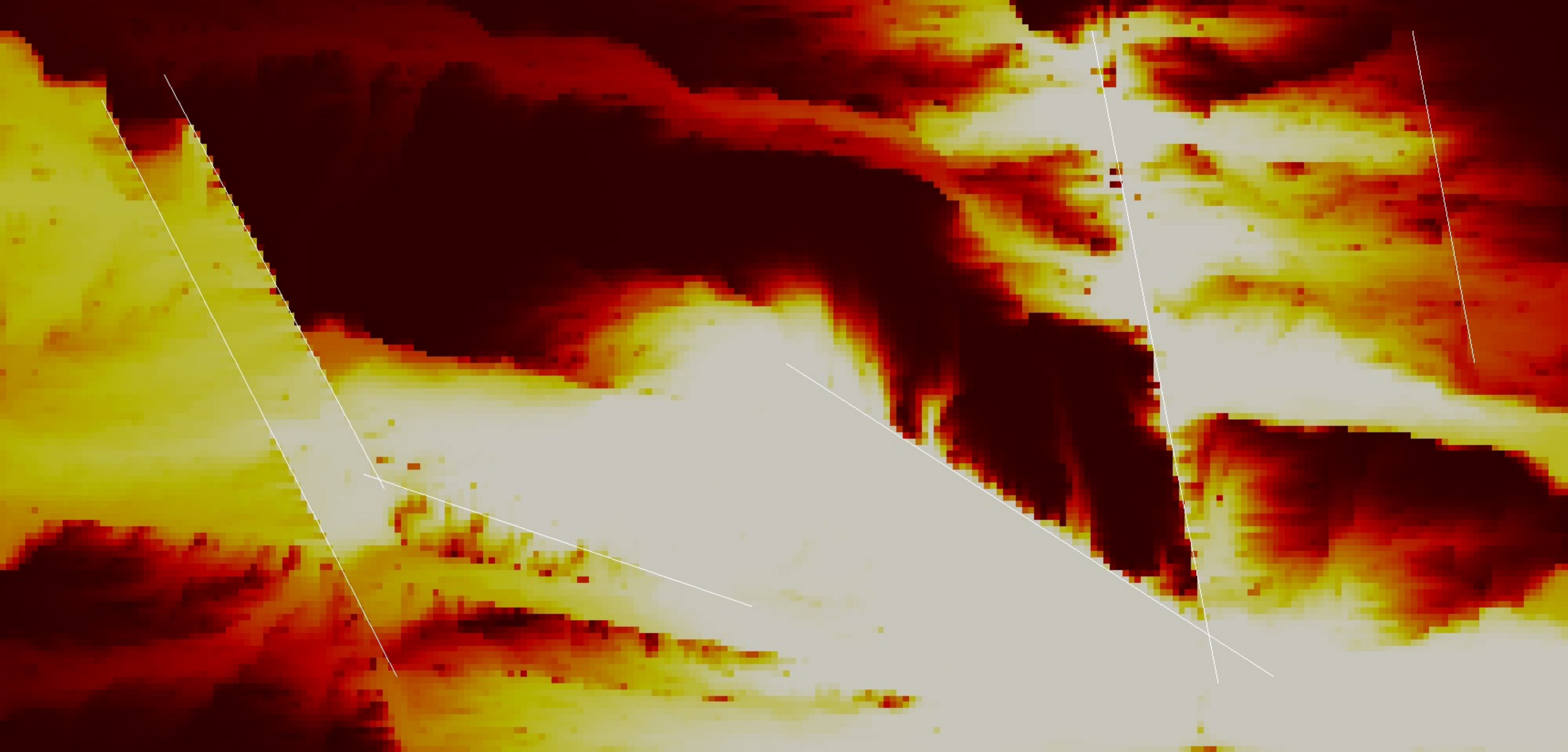}};
		\node[anchor=west] at (.32\linewidth, 2mm) {	  
			
			\pgfplotscolorbardrawstandalone[
			colormap/hot2,
			point meta min=0,
			point meta max=1,
			parent axis height/.initial=2.cm,
			colormap access=map,
			colorbar style={
				title=$s^m$ [-],
				scaled ticks=false,
				minor y tick num = 1,
				/pgf/number format/fixed,
				/pgf/number format/fixed zerofill,
				/pgf/number format/precision=1
			}
			]};
		\end{tikzpicture}
		\caption{}
	\end{subfigure}	  
	\hfill\null
	
	\bigskip
	
	\hfill
	\begin{subfigure}[b]{\panelsize\linewidth}	
		\flushleft    
		
		\begin{tikzpicture}
		\node at (0, 0) {\includegraphics[width=.65\textwidth]{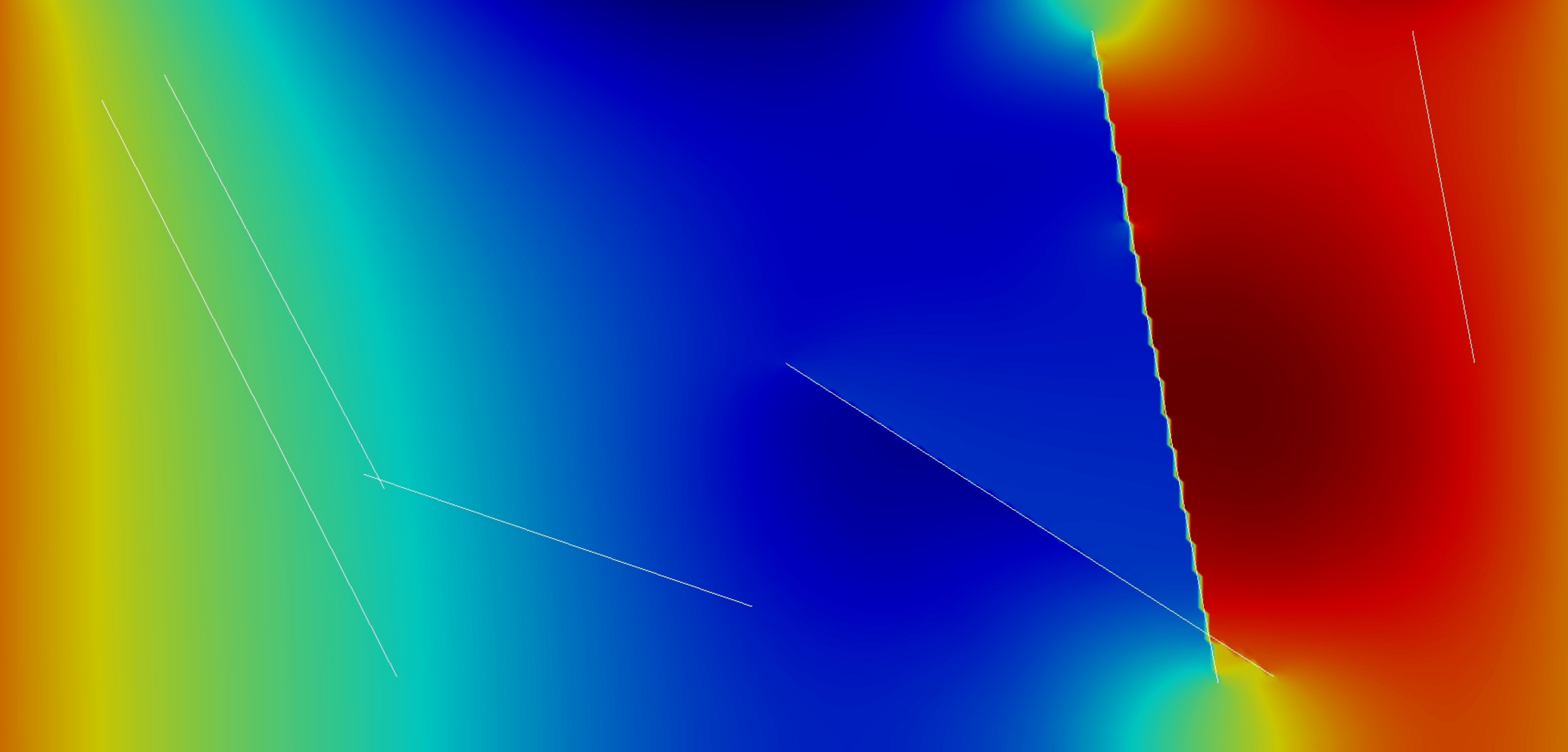}};
		\node[anchor=west] at (.32\linewidth, 2mm) {	     
			
			\pgfplotscolorbardrawstandalone[
			colormap/jet,
			point meta min=-2.752e-1,
			point meta max=1.001e-1,
			parent axis height/.initial=2.cm,
			colormap access=map,
			colorbar style={
				title=$u_x$ [m],
				scaled ticks=false,
				minor y tick num = 3,
				/pgf/number format/sci,
				/pgf/number format/sci zerofill,
				/pgf/number format/precision=1
			}
			]};
		\end{tikzpicture}
		\caption{}
	\end{subfigure}	
	\hfill
	\begin{subfigure}[b]{\panelsize\linewidth}	
		\flushleft     
		
		\begin{tikzpicture}
		\node at (0, 0) {\includegraphics[width=.65\textwidth]{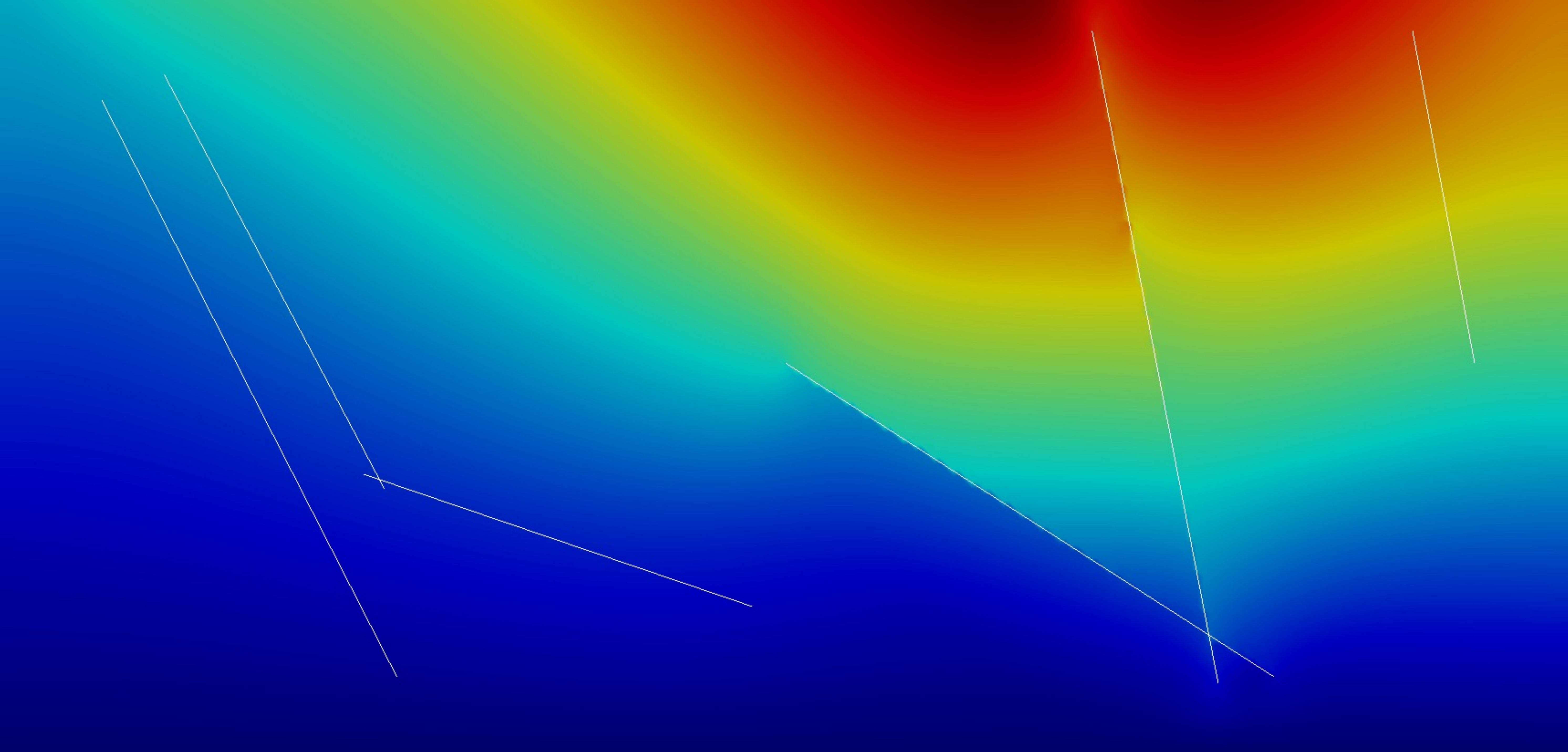}};
		\node[anchor=west] at (.32\linewidth, 2mm) {	      
			
			\pgfplotscolorbardrawstandalone[
			colormap/jet,
			point meta min=0,
			point meta max=4.937e-1,
			parent axis height/.initial=2.cm,
			colormap access=map,
			colorbar style={
				title=$u_y$ [Pa],
				scaled ticks=false,
				minor y tick num = 3,
				/pgf/number format/sci,
				/pgf/number format/sci zerofill,
				/pgf/number format/precision=1
			}
			]};
		\end{tikzpicture}
		\caption{}
	\end{subfigure}	  
	\hfill\null 		
	
	\caption{Test Case 4, Scenario 2: Contour plots of pressure (top left) and saturation (top right), x (bottom left) and y (bottom right) displacements at the end of the simulation.}
	\label{Fig:Case4PandS}
\end{figure}
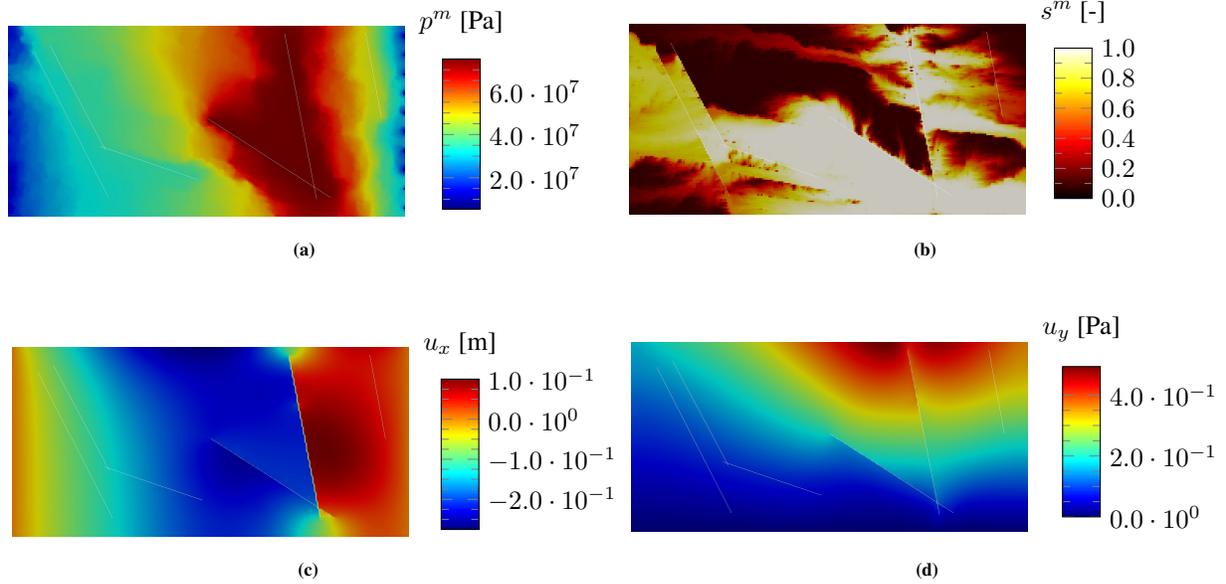

An injector well is now added, injecting water at a constant pressure $p_\text{inj} = 75\,\text{MPa}$. The initial time-step size is 0.01 days and the total simulation time is equal to $T_\text{final} = 1 \, \text{day}$. Figure \ref{Fig:Case4PandS} shows contour plot of pressure, saturation and displacement at the end of the injection process. We observe that the fracture in which water is injected and the one connected to it open significantly due to the pressure increase. Additionally, the saturation path towards the producers is strongly determined by the fracture distribution since all fractures are highly permeable compared to the matrix (even in their closed state). Additionally, the pressure increase creates an uplift at the top boundary which is more pronounced closer to the injection point. 

\section{Concluding remarks}
\label{Sec:Conclusions}
This work has presented an embedded discretization method for coupled multiphase flow and geomechanics in fractured porous media. Fractures are treated as lower dimensional domains embedded in the rock matrix. As such, there is no need for a conforming grid adjusted to the geometry of the fractures.

The linear momentum balance equation is discretized using a finite-element method, whereas a finite-volume method is employed for the mass balance equations.  This mixed method is appealing as it can be readily incorporated into standard reservoir simulation frameworks, and it enforces cell-wise mass conservation. The contribution of fractures to the mechanical and flow problem is captured by employing the EFEM and the EDFM methods, respectively.  Comparisons to an XFEM-based alternative show that all of the methods exhibit similar convergence behavior, with smaller error constants for richer interpolations.

Ongoing research activities are focused on the extension to high-resolution, 3D domains, requiring the development of scalable solver strategies. Future work will also address fracture propagation, with a particular focus on the tip treatment for hydraulic fracturing applications.  

\section*{Acknowledgements}
The authors thank Dr. Chandrasekhar Annavarapu for fruitful discussions. This work was funded by Total S.A. through the FC-MAELSTROM Project.  This work was performed under the auspices of the U.S. Department of Energy by Lawrence Livermore National Laboratory under Contract DE-AC52-07NA27344.

\section*{Appendix: Test Case 4 fracture and well locations}

The end points of each fracture present in test case 4 are presented in table \ref{Tab:FracturesLocation}. Table \ref{Tab:WellsLocation}, instead, summarized the cells perforated by the production and injection wells of test case 4.

\label{AppendixB}
\begin{table}[htbp]
	\centering
	\caption{Test Case 4: Coordinates of the end points of the fractures.}
	\begin{tabular}{c  c}
		\toprule
		\textbf{Fracture} & \textbf{Coordinates} [$(x_1, y_1) - (x_2, y_2)$] \\
		1 & (16.2, 104) -- (63.3, 12)\\
		2 & (26.2, 108) -- (61.3, 42)\\
		3 & (58, 44.3)  -- (120,  23.2)\\
		4 & (174.2, 115) -- (194.3, 11)\\
		5 & (125.3,  62) -- (203.2, 12)\\
		6 & (225.3,  115) -- (235.2, 62)\\
		\bottomrule
	\end{tabular}
	\label{Tab:FracturesLocation}
\end{table}

\begin{table}[htbp]
	\centering
	\caption{Test Case 4: List of cells perforated by each well.}
	\begin{tabular}{c  c}
		\toprule
		\textbf{Well} & \textbf{Perforated cells} \\
		Producer 1 & 1251, 4251, 7251, 10251, 13251, 16251, 19251, 22251, 25251, 28251\\
		Producer 2 & 1500, 4500, 7500, 10500, 13500, 16500, 19500, 22500, 25500, 28500\\
		Injector & Cell 1 of fracture 5\\
		\bottomrule
	\end{tabular}
	\label{Tab:WellsLocation}
\end{table}

\pagebreak

\bibliographystyle{abbrv}
\bibliography{MatteoRef.bib}
\end{document}